\documentclass[smallextended,glov3]{svjour3}
  \def\ZR{{\mathbb R}}

\def\beq{\begin{equation}}
\def\eeq{\end{equation}}
\def\be{\begin{equation}}
\def\ee{\end{equation}}
\def\beqar{\begin{eqnarray}}
\def\eeqar{\end{eqnarray}}
\def\ber{\begin{eqnarray}}
\def\eer{\end{eqnarray}}
\def\berb{\begin{eqnarray*}}
\def\eerb{\end{eqnarray*}}

\def\SO{\mathop{\rm SO}\nolimits}
\def\so{\mathop{\rm so}\nolimits}

\def\det{\mathop{\rm det}\nolimits}

\def\norm#1.#2.{\|#1\|_{#2}}
\def\Norm#1.#2.{\big\|#1\big\|_{#2}}
\def\NOrm#1.#2.{\bigg\|#1\bigg\|_{#2}}
\def\NORm#1.#2.{\Big\|#1\Big\|_{#2}}
\def\NORM#1.#2.{\Bigg\|#1\Bigg\|_{#2}}

\def\vec#1{{\mathchoice{\mbox{\boldmath$\displaystyle#1$}}
{\mbox{\boldmath$\textstyle#1$}}
{\mbox{\boldmath$\scriptstyle#1$}}
{\mbox{\boldmath$\scriptscriptstyle#1$}}}}

\newcommand{\sym}{\mathop{\rm sym}\nolimits}

\def \0b{{\hbox{\boldmath $0$}}}

\newcommand{\ab}{\vec{a}} \newcommand{\bb}{\vec{b}}
\newcommand{\cb}{\vec{c}}

\newcommand{\eb}{\vec{e}} \newcommand{\fb}{\vec{f}}

 \newcommand{\nb}{\vec{n}}
\newcommand{\ovbb}{\vec{o}} 
 \newcommand{\rb}{\vec{r}}
 \newcommand{\tb}{\vec{t}}
\newcommand{\ub}{\vec{u}} \newcommand{\vb}{\vec{v}}
 \newcommand{\xb}{\vec{x}}
\newcommand{\yb}{\vec{y}}

\newcommand{\Abb}{{\bf A}} \newcommand{\Bbbb}{{\bf B}}
\newcommand{\Cbb}{{\bf C}} \newcommand{\Dbb}{{\bf D}}
\newcommand{\Ebb}{{\bf E}} \newcommand{\Fbb}{{\bf F}}
\newcommand{\Gbb}{{\bf G}} 
\newcommand{\Ibb}{{\bf I}} \newcommand{\Jbb}{{\bf J}}
\newcommand{\Kbb}{{\bf K}} \newcommand{\Lbb}{{\bf L}}
\newcommand{\Mbb}{{\bf M}} 
\newcommand{\Obb}{{\bf O}} 
\newcommand{\Qbb}{{\bf Q}} \newcommand{\Rbb}{{\bf R}}

\def \alb{\vec{\alpha}} \def \betab{\vec{\beta}}
\def \gamb{\vec{\gamma}} 
 
 \def \etab{\vec{\eta}}
 
 \def \kab{\vec{\kappa}}

\def \Thetab{\vec{\Theta}}

\def \Thetab{\vec{\Theta}}

\newcommand{\at}{{\tilde a}}

\newcounter{primjer}[section]
\newcounter{tvrdnja}[section]
\newcounter{zadatak}[section]
\newcommand{\oRbb}{{\bf \overline{R}}}

\usepackage{esint}
\usepackage{amssymb}
\usepackage{graphicx}
\usepackage{enumerate}
\newcommand{\dist}{\mathop {\mbox{\rm dist}}\nolimits}

\def \varphib{\vec{\varphi}}

\newcommand{\homeg}{\hat{\Omega}^h}
\newcommand{\wirbb}{\widetilde{\Rbb}}

\def\tra{\mathop{\rm trace}\nolimits}

\begin{document}

\title{Nonlinear weakly curved rod by $\Gamma$-convergence
}




\author{
        Igor Vel\v{c}i\'{c} 
}
\institute{ Igor Vel\v{c}i\'{c} \at Faculty of Electrical Engineering and Computer Science, University of Zagreb, Unska 3, 10000 Zagreb, Croatia \\Tel: +385-1-6129965 \\Fax:+385-1-6170007 \\ \email{igor.velcic@fer.hr}}



\maketitle

\begin{abstract}
We present a nonlinear model of weakly curved rod, namely the
type of curved rod where the curvature is of the order of the
diameter of the cross-section. We use the approach analogous to
the one for rods and curved rods and start from the strain
energy functional of  three dimensional nonlinear elasticity
and do not presuppose any constitutional behavior. To derive
the model, by means of $\Gamma$-convergence,  we need to
propose how is the order of strain energy related to the
thickness of the body $h$. We analyze the situation when the
strain energy (divided by the order of volume) is of the order
$h^4$. That is the same approach as the one when  F\"oppl-von
K\'arm\'an model for plates and the analogous model for rods
are obtained. The obtained model is analogous to Marguerre-von
K\'arm\'an for shallow shells and its linearization is the
linear shallow arch model which can be found in the literature.
\end{abstract}
\keywords{weaky curved rod \and Gamma convergence \and shallow
arch \and asymptotic analysis} \subclass{74K20 \and 74K25}

\section{Introduction}

The study of thin structures is the subject of numerous works
in the theory of elasticity. There is a vast literature on the
subject of rods, plates and shells (see
\cite{Antman,Ciarlet0,Ciarlet1}).

The derivation and justification of the lower dimensional
models, equilibrium and dynamic, of rods, curved rods, weakly
curved rods, plates and shells in linearized elasticity, by
using formal asymptotic expansion, is well established (see
\cite{Ciarlet0,Ciarlet1} and the references therein). In all
these approaches one starts from the equations of
three-dimensional linearized elasticity and then via formal
asymptotic expansion justify the lower dimensional models. One
can also obtain the convergence results. In
\cite{Alvarez1,Alvarez2} the linear model of weakly curved rod
(or as it is called shallow arch) is derived and the
convergence result is obtained. We call weakly curved rods or
shallow arches those characterized by the fact that the
curvature of their centerline should has the same order of
magnitude as the diameter of the cross section, both being much
smaller than their length.

Formal asymptotic expansion is also applied to derive non
linear models of rods, plates and shells (see
\cite{Ciarlet0,Ciarlet1,Marigo} and the references therein),
starting from three-dimensional isotropic elasticity (usually
Saint-Venant-Kirchoff material). Hierarchy of the models is
obtained, depending on the the order of the external loads
related to the thickness of the body $h$ (see also \cite{Fox}
for plates).

However, formal asymptotic expansion does not provide us a
convergence result. The first convergence result, in deriving
lower dimensional models from three-dimensional non linear
elasticity, is obtained applying $\Gamma$-convergence, very
powerful tool introduced by Degiorgi (see
\cite{Braides,dalmaso}). Using $\Gamma$-convergence, elastic
string models, membrane plate and membrane shell models are
obtained (see \cite{Acerbi,Ledret1,LeDret2}). It is assumed
that the external loads are of order $h^0$. The obtained models
are different from those ones obtained by the formal asymptotic
expansion in the sense that additional relaxation of the energy
functional is done.

Recently, hierarchy of models of rods, curved rods, plates and
shells is obtained via $\Gamma$- convergence (see
\cite{Muller0,Muller3,Muller4,Muller6,Lewicka1,Mora0,Mora1,Scardia1,Scardia2}).
Influence of the boundary conditions and the order and  the
type of the external loads is largely discussed for plates (see
\cite{Muller3,Muller5}). Let us mention that
$\Gamma$-convergence results provide us the convergence of the
global minimizers of the total energy functional. Recently,
compensated compactness arguments are used to obtain the
convergence of the stationary points of the energy functional
(see \cite{Mora,Muller7}).

Here we apply the tools developed for rods, plates and shells
to obtain weakly curved rod model by $\Gamma$-convergence. It
is assumed that we have free boundary conditions and that the
strain energy (divided by the order of volume) is of the order
$h^4$, where $h$ is the thickness of the rod. This corresponds
to the situation when external transversal dead loads are of
order $h^3$ (see Remark \ref{napcijela1}). The order $h^4$ of
the strain energy gives F\"oppl-von K\'arm\'an model for
plates, Marguerre-von K\'arm\'an model for shallow shells the
analogous model for rods (see \cite{Muller4,Mora1,Velcic}). The
obtained model is non linear model of the lowest order in the
hierarchy of models and its linearization is shallow arch
model, obtained in \cite{Alvarez1,Alvarez2} for isotropic,
homogenous case (see for comparison Remark \ref{moraremark}
d)). Here we do not presuppose any constitutional behavior and
thus work in a more general framework. The main result is
stated in Theorem \ref{najglavnijiii}.

Throughout the paper  $\bar{A}$ or $\{A\}^-$   denotes the
closure of the set. By a domain we call a bounded open set with
Lipschitz boundary. $\Ibb$  denotes the identity matrix, by
$\SO(3)$ we denote the rotations in $\ZR^3$, by $\so(3)$ the
set of antisymmetric matrices $3 \times 3$ and  $\ZR^{3 \times
3}_{\sym}$ denotes the set of symmetric matrices. By $\sym
\Abb$ we denote the symmetric part of the matrix, $\sym
\Abb=\frac{1}{2}(\Abb+\Abb^T)$. $\eb_1,\eb_2,\eb_3$ are the
vectors of the canonical base in $\ZR^3$. By $\nabla_h$ we
denote $\nabla_h=\nabla_{\eb_1}+\frac{1}{h} \nabla_{\eb_2,
\eb_3}$.  $\| f \|_{C^1(\Omega)}$ stands for $C^1$ norm of the
function $f:\Omega \subset \ZR^n \to\ZR$ i.e. $\|
f\|_{C^1(\Omega)}= \max_{x \in \Omega} |f|+\sum_{i=1}^n\max_{x
\in \Omega} |\partial_i f|$. $\rightarrow$ denotes the strong
convergence and $\rightharpoonup$ the weak convergence.

\section{Setting up the problem}
\setcounter{equation}{0} Let $\omega \subset \ZR^2$ be an open
set
 having area equal to $A$ and
 Lipschitz boundary. For all $h$ such that $0<h \leq 1$ and for
given $L$ we define
\begin{equation}
\omega^h=h \omega, \  \Omega^h=(0,L) \times h \omega.
\end{equation}
We shall leave out superscript when $h=1$, i.e.
$\Omega=\Omega^1$, $\omega=\omega^1$. Let us by $\mu(\omega)$
denote
\begin{equation} \label{o muuu}\mu(\omega)=\int_{\omega} (x_2^2+x_3^2) dx_2
dx_3. \end{equation}
Let us choose coordinate axis such that
\begin{equation} \label{fiksiranje}
\int_{\omega} x_2 dx_2dx_3=\int_{\omega} x_3 dx_2dx_3=\int_{\omega} x_2 x_3 dx_2dx_3=0.
\end{equation}
For every $h$ we define the curve $C^h$ of the form
\begin{equation}
C^h=\{ \theta^h(x_1)=(x_1, \theta^h_2 (x_1), \theta^h_3 (x_1)) \in \ZR^3: x_1 \in (0,L)\}.
\end{equation}
where $\theta^h_k (x_1)$, for $k=2,3$, are given functions
satisfying $\theta^h_k \in C^3(0,L)$. Let $(\tb^h, \nb^h,
\bb^h)$ be the Frenet trihedron associated with the curve $C^h$
\begin{eqnarray}
\tb^h &=& \frac{1}{\sqrt{1+((\theta_2^h)')^2+((\theta_3^h)')^2}}(1,(\theta_2^h)',(\theta_3^h)'), \\
\nb^h &=& \frac{(\tb^h)'}{\|(\tb^h)'\|}, \\
\bb^h &=& \tb^h \times \nb^h.
\end{eqnarray}
We suppose $\nb^h \in C^1 (0,L)$ which is satisfied if
$(\theta^h_2)''$, $(\theta^h_3) ''$ do not vanish at the same
time (which is equivalent to the fact that the curvature of
$C^h$ is strictly positive for any $x_1 \in (0,L)$). The case
where $C^h$ has null curvature points can be treated in the
same fashion, provided that we suppose that along these points
we have the same degree of smoothness as before with $\tb^h$,
$\nb^h$ and $\bb^h$ appropriately chosen (see Remark
\ref{nulazak}). We define the map $\Thetab^h: \bar{\Omega}^h\to
\Thetab^h (\bar{\Omega}^h)=\{\hat{\Omega}^h\}^- \subset \ZR^3$,
where $\hat{\Omega}^h:=\Thetab^h ({\Omega}^h)$,
 in the following manner:
\begin{equation} \label{deftheta}
\Thetab^h (x^h)=(x_1,\theta^h_2(x_1),\theta^h_3(x_1))+x_2^h \nb^h(x_1)+x_3^h \bb^h (x_1)
\end{equation}
and we assume that $\Thetab^h$ is a $C^1$ diffeomorphism which
can be proved if $h$ is small enough and $\theta^h_k$, for
$k=2,3$, are of the form considered here. Namely, we take
$\theta^h_2=h \theta_2$, $\theta^h_3=h \theta_3$ where
$\theta_k \in C^3(0,L)$. Let us suppose \begin{equation}
\label{prepnula} ((\theta_1)'')^2(x_1)+((\theta_2)'')^2(x_1)
\neq 0,
\end{equation}
for all $x_1 \in (0,L)$. A generic point in $\bar{\Omega}^h$ or
$\{\hat{\Omega}^h\}^-$ will be denoted by $x^h=(x_1,
x_2^h,x_3^h)$.

Like in
\cite{Mora0,Mora1,Muller0,Muller3,Muller4,Muller6,Lewicka1} we
start from three dimensional non linear elasticity functional
of strain energy (see \cite{Ciarlet} for an introduction to non
linear elasticity)
\begin{equation} \label{definicijafunkcionala}
I^h (\yb):= \frac{1}{h^2} \int_{\hat{\Omega}^h} W^h(x^h, \nabla \yb)dx^h.
\end{equation}
It is natural to divide the strain energy with $h^2$, since the
volume is vanishing with the order of $h^2$. We are interested
in finding $\Gamma$-limit (in some sense i.e. in characterizing
the limits of minimizers) of the functionals $\frac{1}{h^4}
I^h$. The reason why we divide with $h^4$ is that we want to
obtain theory analogous to F\"oppl-von K\'arm\'an for plates
and rods (see \cite{Muller3,Muller4,Mora1}) and Marguerre-von
K\'arm\'an for shallow shells (see \cite{Velcic}). We do not
look the total energy functional because the part with the
strain energy contains the highest order derivatives (at least
for the external dead loads) and thus makes the most difficult
part of the analysis (see Remarks \ref{napcijela1} and
\ref{napcijela2}). We shall not impose  Dirichlet boundary
condition and assume that the body is free at the boundary. The
consideration of the other boundary conditions is also
possible. We rewrite the functional $I^h$ on the domain
$\Omega$, i.e. we conclude
\begin{equation}
I^h (\yb):= \int_{\Omega} W^h(\Thetab^h\circ P^h (x), (\nabla \yb) \circ \Thetab^h \circ P^h) \det ((\nabla \Thetab^h) \circ P^h(x)) dx,
\end{equation}
where by $P^h:\ZR^3 \to \ZR^3$ we denote the mapping
$P^h(x_1,x_2,x_3)=(x_1,hx_2,hx_3)$.  $(\nabla \yb) \circ
\Thetab^h \circ P^h$ denotes  $\nabla \yb$ evaluated at the
point $\Thetab^h(P^h(x))$. We assume that for each $h$ it is
valid
\begin{equation} \det ((\nabla \Thetab^h) \circ P^h(x)) W^h( \Thetab^h\circ P^h (x), \Fbb) =W(x,\Fbb), \forall x \in \Omega,\  \forall \Fbb \in
\ZR^{3 \times 3}, \end{equation} where the stored energy
function $W$ is independent of $h$ and satisfies the following
assumptions (the same ones as in \cite{Mora1}):
\begin{enumerate}[i)]
\item $W: \Omega \times \ZR^{3 \times 3} \to [0,+\infty]$
    is a Carath\'{e}odory  function; for some $\delta>0$
    the function $\Fbb \mapsto W(x,\Fbb)$ is of class $C^2$
    for $\dist(\Fbb,\SO(3))< \delta$ and for a.e. $x \in
    \Omega$;
\item the second derivative $\frac{\partial^2 W}{\partial
    \Fbb^2}$ is a Carath\'{e}odory function on the set
    $\Omega \times \{ \Fbb \in \ZR^{3 \times 3} :
    \dist(\Fbb, \SO(3)) < \delta \}$ and there exists a
    constant $\gamma>0$ such that
    $$ \Bigg| \frac{\partial^2 W}{\partial \Fbb^2}
    (x,\Fbb)[\Gbb,\Gbb]\Bigg| \leq \gamma |\Gbb|^2 \ \textrm{if } \dist(\Fbb,\SO(3))< \delta \ \textrm{and } \Gbb \in \ZR^{3 \times 3}_{\sym};$$
\item $W$ is frame-indifferent, i.e. $W(x,\Fbb)=W(x, \Rbb
    \Fbb)$ for a.e. $x \in \Omega$ and every $\Fbb \in
    \ZR^{3 \times 3},\Rbb \in \SO(3)$;
\item $W(x,\Fbb)=0$ if $\Fbb \in \SO(3)$; $W(x,\Fbb) \geq C
    \dist^2(\Fbb,\SO(3))$ for every $\Fbb \in \ZR^{3 \times
    3}$, where the constant $C>0$ is independent of $x$.
\end{enumerate}
Under these assumptions we first show the compactness result
(Theorem \ref{prepisano}) i.e. we take the sequence $y^h \in
W^{1,2}(\hat{\Omega}^h;\ZR^3)$ such that
$$ \limsup{h \to 0} \frac{1}{h^4} I^h < +\infty $$
and conclude how that fact affects the limit displacement. In
Lemma \ref{identifikacija}  we prove the lower bound, in
Theorem \ref{upperbound} we prove the upper bound and that
enables us to identify limit functional (Theorem
\ref{najglavnijiii}). First we start with some basic properties
of the mappings $\Thetab^h$ which are necessary for further
analysis.

\section{Properties of the mappings $\Thetab^h$}
\setcounter{equation}{0}

We introduce for $k=2,3$,
\begin{equation}
p_k(x_1) = \frac{\theta_k''(x_1)}{\sqrt{(\theta_2'')^2(x_1)+(\theta_3'')^2(x_1)}}.
\end{equation}
Notice that
\begin{equation} \label{josop}
p_2^2+p_3^2=1, \ p_2p_2'+p_3p_3'=0.
\end{equation}
Let us denote  $p=p_2p_3'-p_2'p_3$.
\begin{theorem} \label{izcea}
Let the functions $\theta^h_k$ be such that
$$\theta^h_k (x_1)=h \theta_k (x_1), \mbox{ {\rm for all} } x_1 \in (0,L), k=2,3.$$
where $\theta_k \in C^3(0,L)$ is independent of $h$. Then there
exists $h_0=h_0(\theta)>0$ such that the Jacobian matrix
$\nabla \Thetab^h(x^h)$, where the mappings $\Thetab^h$ are
defined with (\ref{deftheta}), is invertible for all $x^h \in
\bar{\Omega}^h$ and all $h \leq h_0$. Also there exists $C>0$
such that  for $h \leq h_0$ we have
\begin{equation} \label{determinanta} \det \nabla \Thetab^h=1+ h  \delta^h (x^h), \end{equation}
and
\begin{eqnarray} \label{tangenta}
\tb^h (x_1) &=& \eb_1+h \theta_2 '(x_1) \eb_2+h \theta_3'(x_1) \eb_3 +h^2 \ovbb_1 (x_1), \\ \nonumber
\nb^h (x_1) &=& p_2(x_1) \eb_2+p_3(x_1) \eb_3 -h (\theta_2'p_2+\theta_3' p_3)(x_1) \eb_1\\ \label{normala}& &+h^2 \ovbb_2 (x_1), \\
\nonumber
\bb^h (x_1) &=& -p_3(x_1) \eb_2+p_2(x_1) \eb_3 +h (\theta_2'p_3-\theta_3' p_2)(x_1) \eb_1\\ \label{binormala} & &+h^2 \ovbb_3(x_1),
\\
\nonumber
\nabla \Thetab^h(x^h) &=& \Rbb_e(x_1)+h \Cbb(x_1)+x_2^h \Dbb(x_1)+x_3^h \Ebb (x_1) \\ \label{nablaje} & &+ h^2 \Obb_1^h(x^h),
 \\ \nonumber
(\nabla \Thetab^h(x^h))^{-1} &=& \Rbb_e^T(x_1) -h
\Cbb_1(x_1)-x_2^h \Dbb_1(x_1)-x_3^h \Ebb_1 (x_1) \\
\label{inverz} & &+ h^2 \Obb_2^h(x^h),
\end{eqnarray}
\begin{eqnarray}
\label{glupaocjena1}
\left\| (\nabla \Thetab^h)
-\Rbb_e\right\|_{L^{\infty}(\Omega^h;\ZR^{3 \times 3})} &<& Ch,
\\ \label{glupaocjena2} \left\| (\nabla \Thetab^h)^{-1}
-\Rbb^T_e \right\|_{L^{\infty}(\Omega^h;\ZR^{3 \times 3})}&<&
Ch,
\end{eqnarray}
where
\begin{eqnarray}\label{defR}
\Rbb_e &=& \left( \begin{array}{c|c} 1 & {\bf 0} \\ \hline \\[-2ex] {\bf 0}  & \Rbb
  \end{array} \right), \
\Rbb=\left( \begin{array}{cc} p_2 & - p_3  \\  p_3 & p_2
\end{array} \right), \\ \Cbb&=&\left( \begin{array}{ccc} 0 &
-(\theta_2' p_2+\theta_3' p_3)& -(\theta_3'p_2-\theta_2'p_3)
\\ \theta_2' & 0 & 0 \\ \theta_3' & 0 & 0 \end{array} \right),
\\ \label{defD}
\Dbb&=& \left(\begin{array}{ccc} 0 & 0 & 0 \\ p_2 ' & 0 & 0 \\ p_3 '& 0 & 0   \end{array}  \right), \ \Ebb=\left(\begin{array}{ccc} 0 & 0 & 0 \\ -p_3 ' & 0 & 0 \\ p_2 '& 0 & 0   \end{array}  \right), \\
\Cbb_1&=&\left( \begin{array}{ccc} 0 &   -\theta_2'& -\theta_3'  \\ \theta_2' p_2+\theta_3' p_3 & 0 & 0 \\  \theta_3'p_2-\theta_2'p_3 & 0 & 0 \end{array} \right),\\
 \Dbb_1&=& \left(\begin{array}{ccc} 0 & 0 & 0 \\
0 & 0 & 0 \\  p& 0 & 0   \end{array}  \right), \
\Ebb_1=\left(\begin{array}{ccc} 0 & 0 & 0 \\ -p & 0 & 0 \\
0 & 0 & 0   \end{array}  \right),
\end{eqnarray}
 and $\delta^h: \bar{\Omega}^h \to \ZR$, $\ovbb_i: (0,L) \to
 \ZR^3$, $i=1,2,3$,
$\Obb_k^h:\bar{\Omega}^h \to \ZR^{3 \times 3}$, $k=1,2$ are
functions which satisfy
$$ \sup_{0<h \leq h_0} \max_{x^h \in \bar{\Omega}^h} |\delta^h(x^h)| \leq
C_0,$$  $$\sup_{0<h \leq h_0} \max_{x_1 \in (0,L)}\|\ovbb_i^h(x_1)\| \leq C_0, \ \sup_{0<h \leq h_0} \max_{x_1
\in (0,L)} \|(\ovbb_i^h)'(x_1)\| \leq C_0$$
$$\sup_{0<h \leq h_0} \max_{i,j} \max_{x^h \in \bar{\Omega}^h}
\|\Obb_{k,ij}^h(x^h) \| \leq C_0, \ k=1,2,$$ for some constant
$C_0>0$.
\end{theorem}
\begin{prooof}
It can be easily seen
\begin{equation} \label{tangenta2}
\tb^h(x_1)=\eb_1+h \theta_2 '(x_1) \eb_2+h \theta_3'(x_1) \eb_3 -\frac{h^2}{2}((\theta_2')^2+(\theta_3')^2)\eb_1+h^3\ovbb^h_4(x^h),
\end{equation}
where $\|\ovbb_4\|_{C^2(0,L)} \leq C$. The relations
(\ref{normala}) and (\ref{binormala}) are the direct
consequences of the relation (\ref{tangenta2}). Let us by
$\ub^h:(0,L) \to \ZR^3$ denote the function
\begin{equation}
\ub^h=(1, h \theta_2 ', h \theta_3')^T.
\end{equation}
It is easy to see
\begin{equation} \label{egzaktno}
\nabla \Thetab^h(x^h)=( \ub^h(x_1)+ x_2^h (\nb^h)'(x_1)+x_3^h (\bb^h)' (x_1) \ | \ \nb^h(x_1) \ | \ \bb^h(x_1)).
\end{equation}
The relations (\ref{determinanta}), (\ref{nablaje}),
(\ref{glupaocjena1}), (\ref{glupaocjena2}) are the direct
consequences of the relations (\ref{tangenta}-(\ref{binormala})
and (\ref{egzaktno}). The relation (\ref{inverz}) is the direct
consequence of the fact that, for a regular matrix $\Abb$ and
arbitrary $\Bbbb$, which satisfies $\|\Abb^{-1} \Bbbb\| < 1$
($\|\cdot\|$ is the operational norm), the matrix $\Abb+\Bbbb$
is invertible and
$$\|(\Abb+\Bbbb)^{-1}-(\Abb^{-1}-\Abb^{-1}\Bbbb \Abb^{-1} )\| \leq \frac{\|\Abb^{-1} \Bbbb\|^2 \|\Abb\|^{-1}}{1-\|
\Abb^{-1} \Bbbb\|}. $$
\end{prooof}
To end the proof observe that
\begin{eqnarray*}
\Cbb_1&=& \Rbb^T_e \Cbb \Rbb_e^T, \ \Dbb_1= \Rbb_e^T \Dbb \Rbb_e^T, \  \Ebb_1= \Rbb_e^T\Ebb \Rbb_e^T.
\end{eqnarray*}
\begin{remark} \label{dodatnanapomena}
By a careful computation it can be seen that $\ovbb_2^h$,
$\ovbb_3^h$ and $\ovbb_4^h$ (defined in the relation
(\ref{tangenta2})) are dominantly in $\eb_2, \eb_3$ plane i.e.
that we have for $i=2,3,4$:
\begin{equation}
\ovbb_i^h(x_1)=f_i^2(x_1) \eb_2+f_i^3(x_1) \eb_3+h\rb_i^h(x_1),
\end{equation}
where  $f_i^2,\ f_i^3 \in C^1(0,L)$, $\sup_{0<h \leq
h_0}\|\rb_i^h\|_{C^1(0,L)} \leq C$, for some $C>0$.
\end{remark}
\begin{remark} \label{dodatnanapomena1}
By a further inspection it can be seen that
\begin{eqnarray}
f_4^2&=&-\frac{1}{2} \theta_2'((\theta_2')^2+(\theta_3')^2), \ f_4^3=-\frac{1}{2} \theta_3'((\theta_2')^2+(\theta_3')^2), \\
f_2^2 &=& p_2 \Big(f_2(\theta',\theta'')-\frac{1}{2}((\theta_2')^2+(\theta_3')^2)\Big)-\theta_2'(\theta_2' p_2+\theta_3' p_3), \\ f_2^3 &=& p_3 \Big(f_2(\theta',\theta'')-\frac{1}{2}((\theta_2')^2+(\theta_3')^2)\Big)-\theta_3'(\theta_2' p_2+\theta_3' p_3), \\
f_3^2 &=& p_3((\theta_2')^2+(\theta_3')^2)-p_3f_2(\theta',\theta''),  \\ f_3^3&=&-p_2((\theta_2')^2+(\theta_3')^2)+p_2f_2(\theta',\theta''),
\end{eqnarray}
where $f_2(\theta',\theta'')\in C^1(0,L)$ is the expression
that includes $\theta',\theta''$:
\begin{eqnarray*}
f_2(\theta',\theta'')
&=& \frac{1}{2}\Big((p_2+p_3)((\theta_2')^2+(\theta_3')^2)+2(\theta_2'+\theta_3')(\theta_2'p_2+\theta_3'p_3)\\& &-\sqrt{(\theta_2'')^2+(\theta_3'')^2}(\theta_2'p_2+\theta_3'p_3)^2\Big)
\end{eqnarray*}
\end{remark}
\begin{remark} \label{nulazak}
It is not necessary to impose the condition (\ref{prepnula}).
All we need is the existence of the expansions given by
(\ref{tangenta})-(\ref{binormala}), where $p_2,p_3 \in
C^1(0,L)$, including the statement of Remark
\ref{dodatnanapomena}.
\end{remark}
\begin{remark}
Although $\Thetab^h$ makes the small perturbation of the
central line, \\ $(x_1,0,0)$, for $x_1 \in [0,L]$, it is not
true that $\nabla \Thetab^h$ is close to the identity (like in
the shallow shell model, see \cite{Velcic}). In fact, there are
torsional effects of order 0 on every cross section. This is
the main reason why is the change of coordinates introduced in
the next chapter useful.
\end{remark}

\section{$\Gamma$-convergence}
\setcounter{equation}{0}   We shall need the following theorem
which can be found in \cite{Muller0}.
\begin{theorem}[on geometric rigidity]\label{tgr}
Let $U \subset \ZR^m$ be a bounded Lipschitz domain, $m \geq
2$. Then there exists a constant $C(U)$ with the following
property: for every $\vb\in W^{1,2}(U;\ZR^m)$ there is an
associated rotation $\Rbb\in \SO(m)$ such that
\begin{equation}\label{gr}
\|\nabla \vb-\Rbb\|_{L^2(U)} \leq C(U)\|\dist(\nabla \vb,\SO(m)\|_{L^2(U)}.
\end{equation}
The constant $C(U)$ can be chosen uniformly for a family of domains which are Bilipschitz equivalent with controlled Lipschitz constants.
The constant $C(U)$ is invariant under dilatations.
\end{theorem}
The following version of the Korn's inequality is needed.
\begin{lemma} \label{Kooornova}
Let $\omega \subset \ZR^2$ with Lipschitz boundary and $\ub \in
L^2(\omega;\ZR^2)$. Let us by $e_{ij}(u)$ denote
$e_{ij}(u)=\frac{1}{2}(\partial_i \ub+\partial_j \ub)$. Let us
suppose that for every $i,j=1,2$ we have that $e_{ij}(u) \in
L^2(\omega)$. Then we have that $\ub \in
W^{1,2}(\omega;\ZR^2)$. Also there exists constant $C(\omega)$,
depending only on the domain $\omega$, such that we have
\begin{eqnarray} \nonumber
\|\ub \|_{W^{1,2}(\omega;\ZR^2)} &\leq& C(\omega)\Big(|\int_\omega \ub dx_1 dx_2 |+|\int_\omega (x_1\ub_2-x_2 \ub_1)dx_1dx_2|\\ & &  \label{Koooorn}+\sum_{i,j=1,2}\|e_{ij}(\ub)\|_{L^2(\omega)} \Big).
\end{eqnarray}
Let us suppose that the domains $\omega_s$ are changing in the
sense that they are equal to $\omega_s=\Abb_s \omega$, where
$\Abb_s \in \ZR^{2\times 2}$, and there exists a constant $C$
such that $\|\Abb_s\|,\|\Abb_s^{-1} \| \leq C$. Then the
constant in the inequality (\ref{Koooorn}) can be chosen
independently of $s$.
\end{lemma}
\begin{prooof}
The first part of the lemma (the fixed domain) is a version of
the Korn's inequality (see e.g. \cite{Oleinik}). The last part
we shall prove by a contradiction. Let us suppose the contrary
that for each $n\in \mathbf{N}$ there exists $s^n$ and $\ub^n
\in W^{1,2} (\omega_{s^n};\ZR^2)$ such that we have
\begin{eqnarray}\nonumber
& & |\int_{\omega_{s^n}} \ub^n dx_1 dx_2 |+|\int_{\omega_{s^n}} (x_1\ub^n_2-x_2 \ub^n_1)dx_1dx_2|\\ \label{Koooorn2} & &+\sum_{i,j=1,2}\|e_{ij}^{s^n}(\ub^n)\|_{L^2(\omega_{s^n})} \leq \frac{1}{n} \|\ub^n \|_{W^{1,2}(\omega_{s^n};\ZR^2)},
\end{eqnarray}
where we have by $e_{ij}^{s^n}(\cdot)$ denoted the symmetrized
gradient on the domain $\omega_{s^n}$. Without any loss of
generality we can suppose that $\|\ub^n
\|_{W^{1,2}(\omega_{s^n};\ZR^2)}=1$. Let us take the
subsequence of $(s^n)$ (still denoted by $(s^n)$) such that
$\Abb_{s^n} \to \Abb$ and $\Abb_{s^n}^{-1} \to \Abb^{-1}$ in
$\ZR^{2 \times 2}$.

Let us look the sequence $\ub_c^n=\ub^n \circ \Abb_{s^n} \circ
\Abb^{-1}$. It is clear that there exist $C_1,C_2>0$ such that
\begin{equation} \label{omKorn}
C_1 \leq \| \ub_c^n \|_{W^{1,2}( \omega_\infty ;\ZR^2)} \leq C_2,
\end{equation}
where we have put $\omega_{\infty}:=\Abb \omega$. Thus there
exists $\ub \in W^{1,2} (\omega_\infty ;\ZR^2)$ such that
$\ub_c^n \rightharpoonup \ub$ weakly in $W^{1,2} (\omega_\infty
;\ZR^2)$. Specially, by the compactness of the embedding
$L^2\hookrightarrow W^{1,2}$ (see e.g. \cite{Adams}), we also
conclude the strong convergence $\ub_c^n \to \ub$ in
$L^2(\omega_\infty; \ZR^2)$. Since it is valid $\Abb_{s^n}
\Abb^{-1} \to \Ibb$, it can be easily seen that, from the weak
convergence, it follows $e_{ij}^{s^n}(\ub^n) \circ (\Abb_{s^n}
\circ \Abb^{-1}) \rightharpoonup e_{ij}^{\infty} (\ub)$, weakly
in $L^2 (\omega_\infty;\ZR^2)$, where we have by
$e_{ij}^{\infty}(\cdot)$ denoted the symmetrized gradient on
the domain $\omega_\infty$.  From the weak convergence we can
conclude that
\begin{equation} \label{Koooorn4} \|e_{ij}^{\infty}(\ub)\|_{L^2(\omega_\infty)} \leq
\liminf_{n \to \infty} \|e_{ij}^{s^n}(\ub^n) \circ (\Abb_{s^n}
\circ \Abb^{-1})\|_{L^2(\omega_\infty)}=0, \end{equation} for every
 $i,j=1,2.$ We can also from (\ref{Koooorn2}) conclude that
\begin{equation} \label{Koooorn3} \int_{\omega_\infty} \ub dx_1
dx_2=0, \ \int_{\omega_\infty} (x_1 \ub_2-x_2\ub_1)dx_1 dx_2=0. \end{equation} Applying the standard Korn's inequality on the domain
$\omega_\infty$, i.e.
\begin{eqnarray*}
\|\ub-\ub_c^n\|_{W^{1,2}(\omega_\infty;\ZR^2)} &\leq&
C(\omega_\infty)
\Big(\|\ub-\ub_c^n\|_{L^2(\omega_\infty;\ZR^2)}\\ & &+\sum_{i,j=1,2}
\|e_{ij}(\ub)-e_{ij}(\ub_c^n)\|_{L^2(\omega_{\infty};\ZR^2)}\Big),
\end{eqnarray*}
we conclude that $\ub_c^n \to \ub$ strongly in
$W^{1,2}(\omega_\infty;\ZR^2)$. But then (\ref{omKorn}),
(\ref{Koooorn4}), (\ref{Koooorn3}) make a contradiction with
the version of the Korn's inequality (\ref{Koooorn}) on the
domain $\omega_\infty$.
\end{prooof}

\begin{remark}
The same proof can be done under the assumption that
$\omega_s=F_s(\omega)$, where $F_s$ is the family of
Bilipschitz mappings whose Bilipschitz constants we can control
(i.e. the Lipschitz constants of $F_s$ and $F_s^{-1}$ are
bounded by a universal constant), provided that the family
$F_s$ is strongly compact in $W^{1,\infty}(\omega;\ZR^2)$. It
would require more analysis to conclude the same result only
for Bilipschitz mappings whose Bilipschitz constants we can
control.
\end{remark}

Let us by $\xb':\ZR^3 \to \ZR^3$ denote the change of
coordinates
\begin{equation} \label{zamjenavar}
(x_1', x_2', x_3')= \xb'(x_1,x_2,x_3):= \Rbb_e(x_1) \left( \begin{array}{c} x_1 \\ x_2 \\   x_3 \end{array} \right).
 \end{equation}
  By
 $\Omega'$ we denote $\xb'(\Omega)$ and $\omega'(x_1) \subset \ZR^2$  denotes $\xb'(\{x_1\} \times\omega)$, for $x_1 \in [0,L]$.
The generic point in $\Omega'$ is denoted with
$x'=(x_1',x_2',x_3')$. Let us observe that by (\ref{o muuu})
and (\ref{fiksiranje})
\begin{eqnarray} \label{nulaaa}
& &\int_{\omega'(x_1)} x_2'dx_2'dx_3'=\int_{\omega'(x_1)} x_3'dx_2'dx_3'=0, \\ & & \int_{\omega'(x_1)}x_2'x_3'dx_2'dx_3'=p_2 p_3 \int_\omega(x_2^2-x_3^2)dx_2dx_3,\\ & &
\mu (\omega) =\int_{\omega} (x_2^2+x_3^2) dx_2
dx_3=\int_{\omega'(x_1)} ((x'_2)^2+(x'_3)^2)dx_2' dx_3',
\end{eqnarray}
for all $x_1 \in [0,L]$. By  $(\partial_i \yb_j)\circ \Thetab^h
\circ P^h$)   we  denote  $\partial_i \yb_j$ evaluated at the
point $\Thetab^h(P^h(x))$.

In the sequel we suppose $h_0 \geq 1$ (see Theorem
\ref{izcea}). If this was not the case, what follows could be
easily adapted. Using theorem \ref{tgr} we can prove the
following theorem
\begin{theorem} \label{prepisano}
Let $\yb^h \in W^{1,2}(\hat{\Omega}^h;\ZR^3)$ and let
$$E^h=\frac{1}{h^2} \int_{\homeg}\dist^2 (\nabla \yb^h,\SO(3)) dx.$$
Let us suppose that \begin{equation} \label{limsup} \limsup_{h
\to 0} \frac{E^h}{h^4} <+\infty. \end{equation}

Then there exist maps $\Rbb^h:[0,L] \to \SO(3)$ and
$\widetilde{\Rbb}^h:[0,L] \to \ZR^{3 \times 3}$, with $|\wirbb|
\leq C$, $\wirbb \in W^{1,2}([0,L], \ZR^{3 \times 3})$ and
constants $\oRbb^h \in \SO(3)$, $c^h \in \ZR^3$ such that the
functions $ \widetilde{\yb}^h:= (\oRbb^h)^T\yb^h-c^h$ satisfy
\begin{equation} \label{ocjena000}
\| (\nabla \widetilde{\yb}^h) \circ \Thetab^h\circ P^h -\Rbb^h \|_{L^2(\Omega)} \leq Ch^2,
\end{equation}
\begin{equation} \label{ocjena1}
\|\Rbb^h-\wirbb^h \|_{L^2([0,L])} \leq Ch^2, \quad \|  (\wirbb^h)' \|_{L^2([0,L])} \leq Ch,
\end{equation}
\begin{equation} \label{blizid}
\|\Rbb^h-\Ibb\|_{L^{\infty}([0,L])} \leq Ch.
\end{equation}
Moreover if we define
\begin{eqnarray} \label{defuh}
u^h &=& \frac{1}{A}\int_\omega \frac{\widetilde{\yb}_1^h\circ \Thetab^h\circ P^h-x_1}{h^2}dx_2 dx_3, \\ \label{defvhk}
v^h_k &=& \frac{1}{A}\int_\omega \frac{\widetilde{\yb}_k^h\circ \Thetab^h\circ P^h-h\theta_k}{h} dx_2 dx_3, \ \textrm{for } k=2,3,\\ \label{defwh}
w^h &=& \frac{1}{A\mu(\omega)}\int_\omega \frac{x_2'(\widetilde{\yb}_3\circ \Thetab^h\circ P^h)-x_3'(\widetilde{\yb}_2\circ \Thetab^h\circ P^h)}{h^2} dx_2 dx_3
\end{eqnarray}
then, up to subsequences, the following properties are
satisfied
\begin{enumerate}[(a)]
\item $u^h \rightharpoonup u \ \textrm{in } W^{1,2}(0,L)$;
\item $v_k^h \to v_k \ \textrm{in }  W^{1,2}(0,L)$, where
    $v_k \in W^{2,2}(0,L)$ for $k=2,3$.
\item $w^h \rightharpoonup w \ \textrm{weakly in }
    W^{1,2}(0,L)$;
\item $\frac{(\nabla \widetilde{\yb}^h)\circ \Thetab^h
    \circ P^h-\Ibb}{h} \to \Abb, \textrm{in } L^2(\Omega)$,
    where $\Abb \in W^{1,2} (0,L)$ is given by
\begin{equation} \label{defAbb}
\Abb = \left(\begin{array}{ccc} 0 & -v_2' & -v_3'\\ v_2' & 0& -w \\ v_3' &  w &0\end{array} \right).
\end{equation}

\item $\sym\frac{\Rbb^h-\Ibb}{h^2} \to \frac{\Abb^2}{2} $
    uniformly on $(0,L)$;
\item the sequence $\gamb^h$ defined by
\begin{eqnarray*} \nonumber
\gamb_1^h (x)&=& \frac{1}{h} \Bigg( \frac{(\widetilde{\yb}^h_1\circ \Thetab^h\circ P^h)(x)-x_1}{h^2}-u^h(x_1)\\ & &\hspace{5ex}+x_2'((v_2^h)'+\theta_2')(x_1)+x_3'((v_3^h)'+\theta_3')(x_1) \Bigg),\\
 \nonumber \gamb_k^h (x) &=& \frac{1}{h^2} \Bigg( \frac{(\widetilde{\yb}^h_k\circ \Thetab^h\circ P^h)(x)-h\theta_k-hx_k'}{h}\\ & & \hspace{5ex}-v_k^h(x_1)-h(x'_k)^{\bot}\omega^h(x_1) \Bigg), \ \textrm{for} \ k=2,3,
\end{eqnarray*}
where $(x')^{\bot}:=(0,-x_3',x_2')$, is weakly convergent
in $L^2(\Omega)$ to a function $\gamb$ belonging to the
space $\mathcal{C}$, where
\begin{eqnarray} \nonumber
\mathcal{C}&=&\{ \gamb \in L^2(\Omega;\ZR^3): \int_\omega \gamb=0, \ \partial_2 \gamb, \partial_3 \gamb \in L^2(\Omega; \ZR^3), \\
& & \hspace{-5ex}  \int_{\omega} (x_3' \gamb_2(x_1,\cdot)-x_2' \gamb_3(x_1,\cdot))dx_2 dx_3=0, \ \textrm{for a.e. } x_1 \in (0,L) \}.
\end{eqnarray}
Moreover $\partial_k \gamb^h \rightharpoonup
\partial_k\gamb$ in $L^2 (\Omega)$ for $k=2,3$,
\end{enumerate}
\end{theorem}
\begin{prooof}
We follow the proof of Theorem 2.2 in \cite{Mora1}. Applying
Theorem \ref{tgr} as in the compactness result of \cite{Mora0}
(using the boundedness of $\nabla \Thetab^h$ and $(\nabla
\Thetab^h)^{-1}$  we can find a sequence of piecewise constant
maps $\Rbb^h:[0,L] \to \SO(3)$ such that
\begin{equation}
\int_{\Omega} \| (\nabla \yb^h) \circ \Thetab^h \circ P^h- \Rbb^h \|^2 dx \leq Ch^4,
\end{equation}
and
\begin{equation} \label{josor}
\int_{I'}  \| \Rbb^h (x_1+\xi)-\Rbb(x_1)\|^2 dx_1 \leq Ch^2 (|\xi|+h)^2,
\end{equation}
where $I'$ is any open interval in $(0,L)$ and $\xi \in \ZR$
satisfies $|\xi| \leq \dist(I',\{0,L\})$. Let $\eta \in
C^\infty_0 (0,1)$ be such that $\eta \geq 0$ and $\int_0^1
\eta(s) ds=1$. We set $\eta_h = \frac{1}{h} \eta(\frac{s}{h})$
and we define
$$ \wirbb^h (x_1):= \int_{-h}^h \eta_h (s) \Rbb^h (x_1-s)ds, $$
where we have extended $\Rbb^h$ outside $[0,L]$ by taking
$\Rbb^h (x_1)=\Rbb^h (0)$ for every $x_1<0$, $\Rbb^h(x_1)=
\Rbb^h (L)$ for every $x_1>L$. Clearly $\|\wirbb^h \| \leq C$
for every $h$ while properties (\ref{ocjena1}) follow from
properties (\ref{josor}). Moreover since by construction (see
\cite{Mora0})
$$ \| \Rbb^h (x_1+s)-\Rbb^h (x_1) \|^2 \leq \frac{C}{h^3} \int_{\hat{\Omega}^h}  \dist^2(\nabla \yb^h, \SO(3)) \leq Ch^3, $$
for every $|s| \leq h$ we have by Jensen's inequality that
\begin{equation} \label{ocjbes}
 \| \wirbb^h-\Rbb^h\|^2_{L^\infty([0,L];\ZR^{3 \times 3}} \leq Ch^3
\end{equation}
By the Sobolev-Poincare inequality and the second inequality in
(\ref{ocjena1}), there exist constants $\Qbb^h \in \ZR^{3
\times 3}$ such that $\| \wirbb^h-\Qbb^h \|_{L^{\infty}([0,L];
\ZR^{3 \times 3})} \leq Ch$. Combining this inequality with
(\ref{ocjbes}), we have that $\| \Rbb^h-\Qbb^h
\|_{L^{\infty}([0,L]; \ZR^{3 \times 3})} \leq Ch$. This implies
that $\dist( \Qbb^h, \SO(3)) \leq Ch$; thus, we may assume that
$\Qbb^h$ belongs to $\SO(3)$ and by modifying $\Qbb^h$ by order
$h$, if needed. Now choosing $\oRbb^h =\Qbb^h$ and replacing
$\Rbb^h$ by $(\Qbb^h)^T \Rbb^h$ and $\wirbb^h$ by $(\Qbb^h)^T
\wirbb^h$, we obtain (\ref{blizid}). By suitable choice of
constants $\cb^h \in \ZR^3$ we may assume that
\begin{equation} \label{normal}
\int_{\Omega} (\widetilde{\yb}^h_1 \circ \Thetab^h \circ P^h-x_1)=0, \quad \int_{\Omega} (\widetilde{\yb}^h_k\circ \Thetab^h \circ P^h-h\theta_k) =0, \ \textrm{for } k=2,3.
\end{equation}
Let $\Abb^h =\frac{\Rbb^h-\Ibb}{h}$. By (\ref{blizid}) there
exists $\Abb \in L^{\infty}((0,L);\ZR^{3 \times 3})$ such that,
up to subsequences,
\begin{equation}
\Abb^h \rightharpoonup \Abb \ \textrm{weakly * in } L^{\infty}((0,L);\ZR^{3 \times 3}).
\end{equation}
On the other hand it follows from (\ref{ocjena1}) and
(\ref{blizid}) that
\begin{equation}
\frac{\wirbb^h-\Ibb}{h} \rightharpoonup \Abb \ \textrm{weakly in } W^{1,2} ((0,L);\ZR^{3 \times 3}).
\end{equation}
In particular, $ \Abb \in  W^{1,2} ((0,L);\ZR^{3 \times 3})$
and $h^{-1}(\wirbb^h-\Ibb)$ also converges uniformly. Using
(\ref{ocjbes}) we deduce that
\begin{equation} \label{uniformnoa}
\Abb^h \to \Abb \textrm{ uniformly}.
\end{equation}
In view of (\ref{ocjena000}) this clearly implies the
convergence property in (d). Since $\Rbb^h\in \SO(3)$ we have
$$ \Abb^h+(\Abb^h)^T=-h \Abb^h (\Abb^h)^T. $$
Hence, $\Abb+\Abb^T=0$. Moreover, after division by $2h$ we
obtain property (e) by (\ref{uniformnoa}). For adapting the
proof to the proof of Theorem 2.2~in \cite{Mora1}  it is
essential to see
\begin{eqnarray} \nonumber
(\nabla \widetilde{\yb}^h) \circ \Thetab^h \circ P^h&=&(\nabla(\widetilde{\yb}^h \circ \Thetab^h)\circ P^h) ((\nabla \Thetab^h)^{-1} \circ P^h)\\ \label{kljucnozadokaz1}&=&\nabla_h (\widetilde{\yb}^h \circ \Thetab^h \circ P^h) ((\nabla \Thetab^h)^{-1} \circ P^h).
\end{eqnarray}
From (\ref{kljucnozadokaz1}) it follows
\begin{equation} \label{kljucnozadokaz2}
((\nabla \widetilde{\yb}^h) \circ \Thetab^h \circ P^h)((\nabla \Thetab^h) \circ P^h)=\nabla_h (\widetilde{\yb}^h \circ \Thetab^h \circ P^h).
\end{equation}
and
\begin{eqnarray} \nonumber
(\nabla \widetilde{\yb}^h) \circ \Thetab^h \circ P^h-\Ibb&=&(\nabla_h (\widetilde{\yb}^h \circ \Thetab^h \circ P^h-\Thetab^h \circ P^h))((\nabla \Thetab^h)^{-1} \circ P^h-\Rbb_e^T)\\ \label{zaKorna} & &+(\nabla_h (\widetilde{\yb}^h \circ \Thetab^h \circ P^h-\Thetab^h \circ P^h))\Rbb_e^T.
\end{eqnarray}
Let us notice that from (\ref{deftheta}), (\ref{normala}),
(\ref{binormala}) we can conclude
\begin{equation}
 \label{rastavthetak} \Thetab_k = h\theta_k+hx_k'+O_k(h^3)
\textrm{ for } k=2,3,
\end{equation}
 where
$\|O_k(h^3)\|_{C^1(\Omega)} \leq Ch^3$.

By multiplying (d) with $(\nabla \Thetab^h) \circ P^h=\nabla_h
(\Thetab^h \circ P^h)$ and using (\ref{glupaocjena1}),
(\ref{kljucnozadokaz2}) we obtain
\begin{equation} \label{prilll1}
\frac{\nabla_h (\widetilde{\yb}^h \circ \Thetab^h \circ P^h-\Thetab^h\circ P^h)}{h} \to \Abb \Rbb_e \ \textrm{in } L^2(\Omega).
\end{equation}
Property (b) immediately from (\ref{prilll1}) by using
(\ref{nablaje}), (\ref{normal}) and (\ref{rastavthetak}).
Moreover, $v_k'=\Abb_{k1}$ for $k=2,3$ so that $v_k \in
W^{2,2}(0,L)$ since $\Abb \in W^{1,2} (0,L)$. By using (e),
(\ref{glupaocjena2}), (\ref{ocjena000}) and (\ref{prilll1})
from (\ref{zaKorna}) we conclude that
\begin{equation} \label{zadobu}
\left\| \frac{1}{h^2}\sym (\nabla_h (\widetilde{\yb}^h \circ \Thetab^h \circ P^h-\Thetab^h \circ P^h)\Rbb_e^T) \right\|_{L^2((0,L);\ZR^{3 \times 3})} \leq C
\end{equation}
The weak convergence of $u^h$ follows from (\ref{nablaje}),
(\ref{zadobu}) and
 the definition of $\Rbb_e$. By using the
convergence (\ref{prilll1}) and  Poincare inequality on each
cut $\{x_1\} \times \omega$  we can conclude
\begin{eqnarray} \nonumber
& & \frac{\widetilde{\yb}^h_2 \circ \Thetab^h \circ P^h-\Thetab^h_2 \circ P^h}{h^2}-\frac{1}{h^2A} \int_\omega (\widetilde{\yb}^h_2 \circ \Thetab^h_2 \circ P^h-\Thetab^h_2 \circ P^h) \\ \label{oomegaa1}& & \hspace{10ex} \to (\Abb \Rbb_e)_{22}x_2+(\Abb \Rbb_e)_{23}x_3 \ \textrm{in } L^2(\Omega).
\end{eqnarray}
By using (\ref{fiksiranje}),  (\ref{zamjenavar}) and
(\ref{rastavthetak})
 we
conclude from (\ref{oomegaa1})
\begin{equation} \label{konvergencija111}
w_2^h:=\frac{\frac{1}{h}\widetilde{\yb}^h_2 \circ \Thetab^h \circ P^h-x_2'}{h}-\frac{1}{h^2A} \int_\omega \widetilde{\yb}^h_2 \circ \Thetab^h \circ P^h  \to \Abb_{23} x_3'\ \textrm{in } L^2(\Omega).
\end{equation}

 Let us note that since the left hand side of (\ref{oomegaa1})
i.e. (\ref{konvergencija111}) is bounded in $W^{1,2}(\Omega)$
the convergence in (\ref{konvergencija111}) is in fact weak in
$W^{1,2}(\Omega)$. The only nontrivial thing to prove is the
boundedness of $\partial_1 w_2^h$ in $L^2(\Omega)$. By the
chain rule we have for $i=1,2,3$
\begin{eqnarray}\nonumber
\partial_1(\widetilde{\yb}^h_i \circ \Thetab^h \circ P^h)&=&((\partial_1 \widetilde{\yb}^h_i) \circ \Thetab^h \circ P^h)((\partial_1 \Thetab^h_1) \circ P^h)\\ \label{chainrule1}& &\hspace{-26ex}+((\partial_2 \widetilde{\yb}^h_i) \circ \Thetab^h \circ P^h)((\partial_1 \Thetab^h_2) \circ P^h)+((\partial_3 \widetilde{\yb}^h_i) \circ \Thetab^h \circ P^h)((\partial_1 \Thetab^h_3) \circ P^h)
\end{eqnarray}
and for $k=2,3$

\begin{eqnarray}\nonumber
\partial_k(\widetilde{\yb}^h_i \circ \Thetab^h \circ P^h)&=&h\Big[((\partial_1 \widetilde{\yb}^h_i) \circ \Thetab^h \circ P^h)((\partial_k \Thetab^h_1) \circ P^h)\\ \label{chainrule2}& &\hspace{-27ex}+((\partial_2 \widetilde{\yb}^h_i) \circ \Thetab^h \circ P^h)((\partial_k \Thetab^h_2) \circ P^h)+((\partial_3 \widetilde{\yb}^h_i) \circ \Thetab^h \circ P^h)((\partial_k \Thetab^h_3) \circ P^h)\Big]
\end{eqnarray}

 From (\ref{ocjena000}), (\ref{konvergencija111}) and
(\ref{chainrule1}) we conclude that the boundedness of
$\partial_1 w_2^h$ in $L^2(\Omega)$ is equivalent to the
boundedness of
\begin{eqnarray} \nonumber
z_2^h&=&\frac{\Rbb_{21}^h\partial_1 \Thetab_1 +\Rbb_{22}^h\partial_1 \Thetab_2+\Rbb_{23}^h\partial_1 \Thetab_3-h(p_2'x_2-p_3'x_3)}{h^2}\\ & &-\frac{1}{h^2A} (\Rbb_{21}^h\int_\omega \partial_1 \Thetab_1+\Rbb_{22}^h\int_\omega \partial_1 \Thetab_2+\Rbb_{23}^h\int_\omega \partial_1 \Thetab_3),
\end{eqnarray}
in $L^2(\Omega)$. By using (\ref{fiksiranje}) and
(\ref{egzaktno}) we conclude
\begin{eqnarray} \nonumber
z_2^h &=&\frac{\Rbb_{21}^h(x_2(\nb^h_1)'+x_3(\bb^h_1)') +\Rbb_{22}^h(x_2(\nb^h_2)'+x_3(\bb^h_2)')}{h} \\ \label{josozzz}& &+\frac{\Rbb_{23}^h(x_2(\nb^h_3)'+x_3(\bb^h_3)')-(p_2'x_2-p_3'x_3)}{h}.
\end{eqnarray}
The boundedness of $z_2^h$ in $L^2 (\Omega)$ is the consequence
of (\ref{normala}), (\ref{binormala}) and (\ref{blizid}). Now
we have proved $w_2^h \rightharpoonup \Abb_{23}x'_3$ weakly in
$W^{1,2}(\Omega)$.

Analogously we conclude
\begin{equation}
w^h_3:=\frac{\frac{1}{h}\widetilde{\yb}^h_3 \circ \Thetab^h \circ P^h-x_3'}{h}-\frac{1}{h^2A}\int_\omega \widetilde{\yb}^h_3 \circ \Thetab^h \circ P^h \rightharpoonup -\Abb_{23} x_2',.
\end{equation}
weakly in  $W^{1,2}(\Omega)$.  Now, since $w^h$ can be written
as
\begin{equation}
w^h(x_1)= \frac{1}{A\mu(\omega)} \int_{\omega} (x_2'w^h_3-x_3'w^h_2)dx_2 dx_3,
\end{equation}
it is clear that $w^h$ converges weakly to the function
$w=-\Abb_{23}=\Abb_{32}$ in $W^{1,2} (0,L)$. Let us define for
$\betab^h: \Omega' \to \ZR^3$, $\betab^h=\gamb \circ
\xb'^{-1}$. By the chain rule we have
\begin{eqnarray} \nonumber
\partial_1
\betab^h_i &=& (\partial_1 \gamb^h_i) \circ (\xb')^{-1}+(p_2'x_2'+p_3'x_3')(\partial_2 \gamb^h_i) \circ (\xb')^{-1}\\ \nonumber & &+(-p_3'x_2'+p_2'x_3')(\partial_3 \gamb^h_i) \circ (\xb')^{-1}, \\ \nonumber
\partial_2 \betab^h_i &=&  p_2 (\partial_2 \gamb^h_i) \circ (\xb')^{-1}-p_3  (\partial_3 \gamb^h_i) \circ (\xb')^{-1}, \\ \label{chainrule3}
\partial_3 \betab^h_i &=&   p_3 (\partial_2 \gamb^h_i) \circ (\xb')^{-1} + p_2(\partial_3 \gamb^h_i) \circ (\xb')^{-1}.
\end{eqnarray}
By differentiating $\betab_1$ with respect to $x_k'$, with
k=2,3, we have
\begin{eqnarray}
\partial_2 \betab_1^h &=& \frac{1}{h^3} \partial_2 (\widetilde{\yb}^h_1 \circ \Thetab^h \circ P^h\circ (\xb')^{-1})+\frac{1}{h}((v_2^h)'+\theta_2'), \\
\partial_3 \betab_1^h &=& \frac{1}{h^3} \partial_3 (\widetilde{\yb}^h_1 \circ \Thetab^h \circ P^h \circ (\xb')^{-1})+\frac{1}{h}((v_3^h)'+\theta_3')).
\end{eqnarray}
Let us analyze only $\partial_2 \betab_1^h$. We have by
(\ref{egzaktno}), (\ref{chainrule2}) and the chain rule
\begin{eqnarray} \nonumber
\partial_2 \betab_{1}^h &=& \frac{((\partial_1\widetilde{\yb}^h_1) \circ \Thetab^h \circ P^h\circ (\xb')^{-1}) (p_2\nb^h_1-p_3 \bb^h_1)}{h^2}\\& & \nonumber+\frac{((\partial_2\widetilde{\yb}^h_1) \circ \Thetab^h \circ P^h \circ (\xb')^{-1})(p_2\nb^h_2-p_3 \bb^h_2)}{h^2}\\& & \nonumber+\frac{((\partial_3\widetilde{\yb}^h_1) \circ \Thetab^h \circ P^h\circ (\xb')^{-1}) (p_2\nb^h_3-p_3 \bb^h_3)}{h^2}\\ & &\label{kooo1}+\frac{1}{h}((v_2^h)'+\theta_2').
\end{eqnarray}
By using (\ref{normala}), (\ref{binormala}), (\ref{egzaktno}),
(\ref{ocjena000}), (\ref{blizid}), (\ref{chainrule1}) and the
definition of $v_k^h$ we can conclude that for proving the
boundedness of $\partial_2 \betab_{1}^h$ it is enough to prove
the boundedness of $\delta_{1,2}^h$ in $L^2(\Omega)$ where
\begin{equation} \label{kooo2}
\delta_{1,2}^h:=\frac{-h\Rbb^h_{11}\theta_2'+\Rbb^h_{12}}{h^2}+\frac{\Rbb^h_{21}+h\theta_2'}{h^2}=\frac{1-\Rbb^h_{11}}{h} \theta_2'+\frac{\Rbb^h_{12}+\Rbb^h_{21}}{h^2}.
\end{equation}
 The boundedness
of $\delta_{1,2}^h$ in $L^\infty(\Omega)$ is then the
consequence of the property (e). In the same way we can prove
the boundedness of $\partial_3 \betab_{1}^h$. Using the
Poincare inequality and the fact that $\int_{\omega'(x_1)}
\betab^h_1 dx_2, dx_3'=0$, we deduce that there exists a
constant $C>0$ such that
$$ \int_{\omega'(x_1)} (\betab_1^h(x))^2 dx_2 dx_3 \leq C \int_{\omega'(x_1)}
[(\partial_2 \betab_1^h(x))^2+(\partial_3 \betab_1^h(x))^2]dx_2
dx_3 $$ for a.e. $x_1 \in (0,L)$ and for every $h$. Although
the constant $C$ depends on the domain, since all domains are
translations and rotations of the domain $\omega$, the constant
$C$ can be chosen uniformly. Integrating both sides with
respect to $x_1$, we obtain that the sequence $(\betab_1^h)$ is
bounded in $L^2(\Omega')$ so, up to subsequences $\betab_1^h
\rightharpoonup \betab_1$ and $\partial_k \betab_1^h
\rightharpoonup
\partial_k \betab_1$ weakly in $L^2(\Omega')$, for $k=2,3$.
From the relations (\ref{chainrule3}) it can be concluded that
$\gamb_1^h \rightharpoonup \gamb_1$ and $\partial_k \gamb_1^h
\rightharpoonup
\partial_k \gamb_1$ weakly in $L^2(\Omega)$, for $k=2,3$, where $\gamb = \betab \circ \xb'$.
For the sequences $(\betab_2^h)$, $(\betab_3^h)$, we have by
differentiation that for $j,k=2,3$
\begin{equation}
\partial_j \betab_k^h =\frac{1}{h^2} \Bigg( \frac{1}{h}\partial_j(\widetilde{\yb}^h_k \circ \Thetab^h \circ P^h\circ (\xb')^{-1})-h\delta_{jk} -hw^h(1-\delta_{jk})(-1)^k\Bigg).
\end{equation}
By using the chain rule we see that for $k=2,3$,
\begin{eqnarray*}
\partial_2(\widetilde{\yb}^h_k \circ \Thetab^h \circ P^h\circ (\xb')^{-1})&=&h\Big(((\partial_1\widetilde{\yb}^h_k) \circ \Thetab^h \circ P^h\circ (\xb')^{-1}) (p_2\nb^h_1-p_3 \bb^h_1)\\ &+& ((\partial_2\widetilde{\yb}^h_k) \circ \Thetab^h \circ P^h\circ (\xb')^{-1})(p_2\nb^h_2-p_3 \bb^h_2) \\&+& ((\partial_3\widetilde{\yb}^h_k) \circ \Thetab^h \circ P^h\circ (\xb')^{-1})(p_2\nb^h_3-p_3 \bb^h_3)\Big),\\
\partial_3(\widetilde{\yb}^h_k \circ \Thetab^h \circ P^h\circ (\xb')^{-1})&=& h\Big(((\partial_1\widetilde{\yb}^h_k) \circ \Thetab^h \circ P^h\circ (\xb')^{-1}) (p_3\nb^h_1+p_2 \bb^h_1)\\ &+&((\partial_2\widetilde{\yb}^h_k) \circ \Thetab^h \circ P^h\circ (\xb')^{-1})(p_3\nb^h_2+p_2 \bb^h_2) \\&+&((\partial_3\widetilde{\yb}^h_k) \circ \Thetab^h \circ P^h\circ (\xb')^{-1})(p_3\nb^h_3+p_2 \bb^h_3)\Big).
\end{eqnarray*}
Now we want to check that for $j,k=2,3$
\begin{equation}
e_{jk} (\betab^h):= \frac{1}{2} (\partial_j \betab_k^h+\partial_k \betab_j^h).
\end{equation}
is bounded  in $L^2(\Omega')$. In the similar way as for
$\betab_1$ (relations (\ref{kooo1}) and (\ref{kooo2})) we can
using (\ref{normala}), (\ref{binormala}), (\ref{ocjena000}),
(\ref{blizid}) and the property (e)  conclude that for every
$j,k=2,3$, $e_{jk} (\betab^h) \in L^2(\Omega')$. By using
Korn's inequality (Lemma \ref{Kooornova})  we have that there
exists $C>0$ such that
\begin{eqnarray} \nonumber & &
\| \betab^h_2 \|^2_{W^{1,2}(\omega'(x_1))}+\| \betab^h_3 \|^2_{W^{1,2}(\omega'(x_1))} \leq C \Big( |\int_{\omega'(x_1)}\betab_2^hdx_2'dx_3'|\\ \nonumber & & \hspace{3ex}+|\int_{\omega'(x_1)}\betab_3^hdx_2'dx_3'|+|\int_{\omega'(x_1)} (x_3' \betab_2^h-x_2' \betab_3^h)dx_2'dx_3'|\\ & & \label{Kooorn} \hspace{3ex}+ \sum_{j,k=1,2} \| e_{jk} (\betab^h) \|_{L^2(\omega'(x_1))}\Big),
\end{eqnarray}
for a.e. $x_1 \in (0,L)$.  From the definition of $v_k^h$ and
$w^h$ we see that the functions ($\betab_2^h(x_1,\cdot)$,
$\betab_3^h(x_1,\cdot)$)  belong to the space
\begin{eqnarray} \nonumber
\mathcal{B}_{x_1}&=&\{ \betab=(\betab_2,\betab_3) \in W^{1,2}(\omega'(x_1);\ZR^2): \int_{\omega'(x_1)} \betab dx_2'dx_3'=0, \\ & & \hspace{15ex} \int_{\omega'(x_1)} (x_2'\betab_3-x_3' \betab_2)dx_2' dx_3'=0 \}
\end{eqnarray}
for every $x_1$.
 By integrating
(\ref{Kooorn}) with respect to $x_1$ we conclude that
$\betab_2^h$,$\betab_3^h$ are bounded in $L^2(\Omega')$ as well
as their derivatives with respect to $x_2,x_3$. From this we
can conclude the same fact about $\gamb_2^h$, $\gamb_3^h$. The
fact that the weak limit belongs to the space $\mathcal{C}$ can
be concluded from the fact that for every $h$ and a.e. $x_1$
$(\betab_2^h(x_1,\cdot), \betab_3^h(x_1,\cdot)) \in
\mathcal{B}_{x_1}$. This finishes the proof of (f).

\end{prooof}

\subsection{Lower bound}
\begin{lemma} \label{identifikacija}
Let $\yb^h$, $\widetilde{\yb}^h$, $E^h$, $\Rbb^h$, $u^h$,
$v^h$, $w^h$, $\gamb^h$, $\betab^h=\gamb^h \circ (\xb')^{-1}$,
$\gamb$,  $\betab=\gamb \circ (\xb')^{-1}$, $\Abb$ be as in
Theorem \ref{prepisano} and let us suppose that the condition
(\ref{limsup}) is satisfied and that $\gamb^h \rightharpoonup
\gamb$, $\partial_2 \gamb^h \rightharpoonup
\partial_2 \gamb$, $\partial_3 \gamb^h \rightharpoonup
\partial_3 \gamb$ weakly in $L^2(\Omega)$ i.e.
$\betab^h \rightharpoonup \betab$, $\partial_2 \betab^h
\rightharpoonup
\partial_2 \betab$, $\partial_3 \betab^h \rightharpoonup
\partial_3 \betab$ weakly in $L^2(\Omega')$.
Let us define
\begin{eqnarray} \nonumber
\etab_1^h (x) &=& \frac{1}{h} \Bigg( \frac{(\widetilde{\yb}^h_1\circ \Thetab^h\circ P^h)(x)-\Thetab^h_1 \circ P^h}{h^2}-u^h(x_1)\\ \label{defetah1}& &\hspace{5ex}+x_2'(v_2^h)'(x_1)+x_3'(v_3^h)'(x_1) \Bigg),\\
 \nonumber \etab_k^h (x) &=& \frac{1}{h^2} \Bigg( \frac{(\widetilde{\yb}^h_k\circ \Thetab^h\circ P^h)(x)-\Thetab^h_k \circ P^h}{h}\\ & & \label{defetahk} \hspace{5ex}-v_k^h(x_1)-h(x'_k)^{\bot}\omega^h(x_1) \Bigg), \ \textrm{for} \ k=2,3,
\end{eqnarray}
and $\kab^h=\etab^h \circ (\xb')^{-1}$. Then we have that
$\etab^h \rightharpoonup \etab$ weakly in $L^2(\Omega)$ and
$\partial_k \etab^h \rightharpoonup \partial_k \etab$ weakly in
$L^2(\Omega)$ i.e. $\kab^h \rightharpoonup \kab$, $\partial_2
\kab^h \rightharpoonup
\partial_2 \kab$, $\partial_3 \kab^h \rightharpoonup
\partial_3 \kab$ weakly in $L^2(\Omega')$.
Here
\begin{eqnarray} \label{defeta1}
\etab_1&=&\gamb_1, \\
\etab_2 &=& \gamb_2+f_2^2x_2+f_2^3 x_3=\gamb_2+g_2^2 x_2'+g_2^3 x_3', \\ \label{defeta3}
\etab_3 &=& \gamb_3+f_3^2x_2+f_3^3 x_3=\gamb_3+g_3^2 x_2'+g_3^3 x_3', \\ \label{defkab}
\kab &=&\etab \circ (\xb')^{-1},
\end{eqnarray}
$f_k^j$ are defined in Remark \ref{dodatnanapomena1} and
$g_k^j$ can be easily defined for the above identities to be
valid i.e. for $k=2,3$, we define
\begin{equation}
g_k^2=p_2f_k^2-p_3f_k^3, \ g_k^3=p_3 f_k^2+p_2 f_k^3.
\end{equation}
The following strain convergence is valid
\begin{equation} \label{strainkon}
\Gbb^h:= \frac{(\Rbb^h)^T ((\nabla \yb^h) \circ \Thetab^h\circ P^h)-\Ibb}{h^2} \rightharpoonup \Gbb \quad \textrm{in } L^2(\Omega;\ZR^{3 \times 3}).
\end{equation}
and the symmetric part of $\Gbb$ denoted by $\widetilde{\Gbb}$,
satisfies
\begin{equation}
\widetilde{\Gbb}=\sym(\Jbb-\frac{1}{2} \Abb^2+ \Kbb),
\end{equation}
where
\begin{eqnarray}  \label{defJbb}
\Jbb&=& \left(\begin{array}{ccc} u'+v_2'\theta_2'+v_3'\theta_3'
& 0 & 0 \\ w \theta_3'  & v_2' \theta_2'& v_2'\theta_3' \\ -w
\theta_2'& v_3'\theta_2' & v_3'\theta_3'   \end{array} \right)\\ \label{defKbb}
\Kbb &=& \left(\begin{array}{c}-x_2' v_2''-x_3'v_3'' \\ -x_3'w'\\
x_2'w'\end{array} \ \Bigg| \  \partial_2 \kab \ \Bigg| \  \partial_3 \kab \right).
\end{eqnarray}
Moreover,
\begin{eqnarray*} \liminf_{h \to 0} \frac{1}{h^6} \int_{\hat{\Omega}^h} W^h(x,
\nabla \hat{\yb}^h)dx&=& \liminf_{h \to 0}
\frac{1}{h^4}\int_{\Omega} W(x^h,(\nabla \hat{\yb}^h) \circ
\Thetab^h \circ P^h)dx^h\\&\geq&  \frac{1}{2} \int_{\Omega}
Q_3(x,\widetilde{\Gbb}(x)) dx,\end{eqnarray*} where $Q_3$ is twice the quadratic form of linearized elasticity, i.e.,
\begin{equation} \label{defQ3}
Q_3 (x,\Fbb)= \frac{\partial^2 W}{\partial \Fbb^2} (\Ibb)[\Fbb,\Fbb].
\end{equation}
\end{lemma}
\begin{prooof} We follow the proof of Lemma 2.3 in \cite{Mora1}.
Firstly, using Remark \ref{dodatnanapomena}, it can be seen
that
\begin{eqnarray*}
\etab_1^h&=&\gamb^h_1+h\ovbb_1,\\
\etab_2^h &=& \gamb^h_2+f_2^2x_2+f_2^3x_3+h\ovbb_2, \\
\etab_3^h &=& \gamb^h_3+f_3^2 x_2+f_3^3 x_3+h\ovbb_3,
\end{eqnarray*}
where $\|\ovbb_i\|_{C^1(\Omega)} \leq C$, for some $C>0$. The
convergence of $\etab^h$ is an easy consequence of the
convergence of $\gamb^h$.  The estimate (\ref{ocjena000})
implies that the $L^2$ norm of $\Gbb^h$ is bounded; therefore
up to subsequences, there exists $\Gbb\in L^2(\Omega; \ZR^{3
\times 3})$ such that (\ref{strainkon}) is satisfied. In order
to identify the symmetric part of $\Gbb$ we decompose $\Rbb^h
\Gbb^h$ as follows:
\begin{equation}
\Rbb^h \Gbb^h= \frac{(\nabla \widetilde{\yb}^h) \circ \Thetab^h \circ P^h- \Ibb}{h^2}-\frac{\Rbb^h-\Ibb}{h^2},
\end{equation}
so that
\begin{equation} \label{konF}
\Fbb^h := \sym \frac{(\nabla \widetilde{\yb}^h) \circ \Thetab^h \circ P^h- \Ibb}{h^2}= \sym (\Rbb^h \Gbb^h)+ \sym \frac{\Rbb^h-\Ibb}{h^2}.
\end{equation}
The right hand side converges weakly to
$\widetilde{\Gbb}+\frac{\Abb^2}{2}$ by (\ref{blizid}),
(\ref{strainkon}) and property (e) of the Theorem
\ref{prepisano}. Therefore the sequence $\Fbb^h$ has a weak
limit in $L^2(0,L)$, satisfying
$\Fbb=\widetilde{\Gbb}+\frac{\Abb^2}{2}$. To conclude we only
need to identify $\Fbb$. Consider the functions
\begin{equation}
\phi_1^h := \frac{\widetilde{\yb}_1^h \circ \Thetab^h \circ P^h-x_1}{h^2}.
\end{equation}
From property (f) of Theorem \ref{prepisano} it follows that
the functions
$\phi_1^h-u^h+x_2'((v_2^h)'+\theta_2')+x_3'((v_3^h)'+\theta_3')$,
which are equal to $h\gamma_1^h$ converge strongly to $0$ in
$L^2(\Omega)$. Thus by property (a) and (b) of Theorem
\ref{prepisano} we conclude that
\begin{equation} \label{dosada1}
\phi_1^h \to u-x_2'(v_2'+\theta_2')-x_3'(v_3'+\theta_3') \ \textrm{in } L^2(\Omega).
\end{equation}
 By using the chain rule, the property (d) of
Theorem \ref{prepisano}, (\ref{normala}), (\ref{binormala}),
(\ref{rastavthetak}) we can conclude that
\begin{equation} \label{dosada2}
\partial_1 \phi_1^h \rightharpoonup \Fbb_{11}-x_2(\theta_2'p_2+\theta_3'p_3)'+x_3(\theta_2'p_3-\theta_3'p_2)'-v_2'(\partial_1 x_2'+\theta_2')-v_3'(\partial_1 x_3'+\theta_3'),
\end{equation}
weakly in $L^2(\Omega)$.  From (\ref{dosada1}) and
(\ref{dosada2}) we conclude that
\begin{eqnarray} \nonumber
& &u'-\partial_1 x_2'(v_2'+\theta_2')-x_2'(v_2''+\theta_2'')-\partial_1 x_3'(v_3'+\theta_3')-x_3'(v_3''+\theta_3'') \\ \nonumber &= &
\Fbb_{11}-x_2(\theta_2'p_2+\theta_3'p_3)'+x_3(\theta_2'p_3-\theta_3'p_2)'\\ & &-v_2'(\partial_1 x_2'+\theta_2')-v_3'(\partial_1 x_3'+\theta_3').
\end{eqnarray}
After some calculation we obtain
\begin{equation}
\Fbb_{11}=u'+v_2'\theta_2'+v_3'\theta_3'-x_2'v_2''-x_3'v_3''.
\end{equation}
To identify $\Fbb_{12}$ we have to do some straight forward
computations. By using the chain rule, (\ref{normala}),
(\ref{binormala}), Remark \ref{dodatnanapomena}, property (d)
of Theorem (\ref{prepisano}) we can conclude
\begin{eqnarray} \nonumber
& &\frac{1}{h^2} \partial_1 (\widetilde{\yb}^h_2 \circ \Thetab^h \circ P^h)+\frac{1}{h^3} \Big(p_2\partial_2 (\widetilde{\yb}_1^h \circ \Thetab^h \circ P^h)-p_3\partial_3 (\widetilde{\yb}_1^h \circ \Thetab^h \circ P^h)\Big)= \\ \label{pomiz1} & &
2\Fbb_{12}^h+\frac{\partial_1 x_2'}{h}-w(\theta_3'+\partial_1 x_3')+\Obb_1^h,
\end{eqnarray}
where $\lim_{h \to 0}\|\Obb_1^h\|_{L^2(\Omega;\ZR^{3 \times
3})}=0$. On the other hand it can be easily seen that
\begin{eqnarray} \nonumber
\partial_2 \betab_1^h &=& p_2 \partial_2 \gamma_1^h-p_3 \partial_3 \gamma_1^h\\ \nonumber
&=& \frac{1}{h^3} \Big(p_2\partial_2 (\widetilde{\yb}_1^h \circ \Thetab^h \circ P^h)-p_3\partial_3 (\widetilde{\yb}_1^h \circ \Thetab^h \circ P^h)\Big)\\ & &  \label{pomiz2}+\frac{1}{h}((v_2^h)'+\theta_2').
\end{eqnarray}
From (\ref{konvergencija111}), (\ref{pomiz1}), (\ref{pomiz2})
we conclude
\begin{eqnarray} \nonumber
2\Fbb^h_{12}&=&\frac{1}{h^2} \partial_1 (\widetilde{\yb}^h_2 \circ \Thetab^h \circ P^h-hx_2')-\frac{1}{h}((v_2^h)'+\theta_2') \\ \nonumber& & +\partial_2 \betab_1^h+w\theta_3'+w \partial_1 x_3'-\Obb_1^h \\ \label{imamga}
&=&\partial_1 w_2^h+\partial_2 \betab_1^h+w\theta_3'+w \partial_1 x_3'-\Obb_1^h.
\end{eqnarray}
By using (\ref{konvergencija111}) we conclude that the right
hand side of (\ref{imamga}) converges in $W^{-1,2}(\Omega)$ to
\begin{equation}
\partial_1 (-wx_3')+\partial_2 \betab_1+w\theta_3'+w \partial_1 x_3'=-x_3'w'+w\theta_3'+\partial_2 \kab_1,
\end{equation}
since $\betab_1=\kab_1$. On the other hand we know that the
left hand side of (\ref{imamga}) converges strongly in
$L^2(\Omega)$ to $2\Fbb_{12}$ and thus we can conclude
\begin{equation}
\Fbb_{12}=\frac{1}{2}(-x_3'w'+w\theta_3'+\partial_2 \kab_1).
\end{equation}
In the same way one can prove
\begin{equation}
\Fbb_{13}=\frac{1}{2}(x_2'w'-w\theta_2'+\partial_3 \kab_1).
\end{equation}
To identify $\Fbb_{22}$ let us observe that by the chain rule,
(\ref{normala}), (\ref{binormala}) and the property (d) of
Theorem (\ref{prepisano}) we have
\begin{eqnarray} \nonumber
& &\frac{1}{h^3}\Big(p_2 \partial_2 (\widetilde{\yb}^h_2 \circ \Thetab^h \circ P^h-\Thetab_2 \circ P^h)-p_3 \partial_3 (\widetilde{\yb}^h_2 \circ \Thetab^h \circ P^h-\Thetab_2 \circ P^h)\Big)=\\
& & \label{pomiz3} \Fbb_{22}^h-v_2'\theta_2'+\Obb_2^h,
\end{eqnarray}
where $\lim_{h \to 0} \| \Obb_2^h \|_{L^2(\Omega;\ZR^{3 \times
3})}=0$. On the other hand we can conclude
\begin{eqnarray} \nonumber
& & \partial_2\kab^h_2 = p_2 \partial_2 \etab^h_2- p_3 \partial_3 \etab^h_3 \\ \nonumber & &= \frac{1}{h^3}\Big(p_2 \partial_2 (\widetilde{\yb}^h_2 \circ \Thetab^h \circ P^h-\Thetab_2 \circ P^h)-p_3 \partial_3 (\widetilde{\yb}^h_2 \circ \Thetab^h \circ P^h-\Thetab_2 \circ P^h)\Big). \\
\end{eqnarray}
In the same way as before we conclude that
\begin{equation}
\Fbb_{22}=v_2' \theta_2'+\partial_2 \kab_2.
\end{equation}
Analogously we can conclude
\begin{equation}
\Fbb_{33}=v_3' \theta_3'+\partial_3 \kab_3.
\end{equation}
To identify $\Fbb_{23}=\Fbb_{32}$ we, by using the chain rule,
(\ref{normala}), (\ref{binormala}) and the property (d) of
Theorem (\ref{prepisano}), can conclude:
\begin{eqnarray} \nonumber
& &\frac{1}{h^3}\Big(p_2 \partial_2 (\widetilde{\yb}^h_3 \circ \Thetab^h \circ P^h-\Thetab_2 \circ P^h)-p_3 \partial_2 (\widetilde{\yb}^h_3 \circ \Thetab^h \circ P^h-\Thetab_2 \circ P^h)\Big)= \\ \label{pomiz6}& & \frac{1}{h^2}(\partial_2 \widetilde{\yb}^h_3) \circ \Thetab^h \circ P^h-v_2' \theta_3' +\Obb_3^h,
\end{eqnarray}
where $\lim_{h \to 0} \| \Obb_3^h \|_{L^2(\Omega;\ZR^{3 \times
3})}=0$. In the same way we conclude
\begin{eqnarray} \nonumber
& &\frac{1}{h^3}\Big(p_3 \partial_2 (\widetilde{\yb}^h_2 \circ \Thetab^h \circ P^h-\Thetab_2 \circ P^h)+p_2 \partial_3 (\widetilde{\yb}^h_2 \circ \Thetab^h \circ P^h-\Thetab_2 \circ P^h)\Big)= \\ & & \frac{1}{h^2}(\partial_2 \widetilde{\yb}^h_3) \circ \Thetab^h \circ P^h-v_3' \theta_2' +\Obb_4^h,
\end{eqnarray}
where $\lim_{h \to 0} \| \Obb_4^h \|_{L^2(\Omega;\ZR^{3 \times
3})}=0$. It can be also concluded
\begin{eqnarray} \nonumber
p_2 \partial_2 \etab_3^h-p_3\partial_3 \etab_3^h &=&\frac{1}{h^3}\Big( p_2 \partial_2 (\widetilde{\yb}^h_3 \circ \Thetab^h \circ P^h-\Thetab_2 \circ P^h)\\ & &-p_3 \partial_2 (\widetilde{\yb}^h_3 \circ \Thetab^h \circ P^h-\Thetab_2 \circ P^h)\Big)+\frac{1}{h^2}w^h,\\
\nonumber p_3 \partial_2 \etab_2^h+p_2\partial_3 \etab_2^h &=& \frac{1}{h^3}\Big(p_3 \partial_2 (\widetilde{\yb}^h_2 \circ \Thetab^h \circ P^h-\Thetab_2 \circ P^h)\\ \label{pomiz7} & &+p_2 \partial_2 (\widetilde{\yb}^h_3 \circ \Thetab^h \circ P^h-\Thetab_2 \circ P^h)\Big)-\frac{1}{h^2}w^h.
\end{eqnarray}
By summing the relations (\ref{pomiz6})-(\ref{pomiz7}) and
letting $h \to 0$ it can be concluded that
\begin{equation}
2\Fbb_{23}= v_2' \theta_3'+v_3' \theta_2'+\partial_2 \kab_3+\partial_3 \kab_2.
\end{equation}
To prove the lower bound we can continue in the same way as in
the proof of Lemma 2.3 in \cite{Mora1}, by using the Taylor
expansion, the cutting and Scorza-Dragoni theorem.
\end{prooof}

\subsection{Upper bound}
 \begin{theorem}[optimality of lower bound] \label{upperbound}
Let $u,w \in W^{1,2}(0,L)$ and $v_k \in W^{2,2} (0,L)$  for
$k=2,3$. Let $\gamb$ be a function in $\mathcal{C}$ where
\begin{eqnarray} \nonumber
\mathcal{C}&=&\{ \gamb \in L^2(\Omega;\ZR^3): \int_\omega \gamb=0, \ \partial_2 \gamb, \partial_3 \gamb \in L^2(\Omega; \ZR^3), \\
& &  \int_{\omega} (x_3' \gamb_2(x_1,\cdot)-x_2' \gamb_3(x_1,\cdot))dx_2 dx_3=0, \ \forall x_1 \in (0,L) \}.
\end{eqnarray}
Set
\begin{equation}
\widetilde{\Gbb}=\sym(\Jbb-\frac{1}{2} \Abb^2+ \Kbb).
\end{equation}
Here $\Abb$, $\Jbb$, $\Kbb$ are defined by the expressions
(\ref{defAbb}), (\ref{defJbb}) and (\ref{defKbb}) and $\etab$,
$\kab$ are defined by the expressions
(\ref{defeta1})-(\ref{defkab}).

Then there exists a sequence $(\hat{\yb}^h) \subset W^{1,2}
(\hat{\Omega}^h,\ZR^3)$ such that for $u^h,v^h_k,w$ defined by
the expressions (\ref{defuh})-(\ref{defwh}) the properties
(a)-(d) of Theorem \ref{prepisano} are valid. Also we have that
the property (f) of Theorem \ref{prepisano} is valid (which is
equivalent that for $\etab^h$ defined by the expressions
(\ref{defetah1})-(\ref{defetahk}) it is valid $\etab^h
\rightharpoonup \etab$ weakly in $L^2(\Omega)$ and $\partial_k
\etab^h \rightharpoonup \partial_k \etab$ weakly in
$L^2(\Omega)$). Also the following convergence is valid

\begin{equation} \label{upperjed}
\lim_{h \to 0} \frac{1}{h^6} \int_{\hat{\Omega}^h} W^h(x^h, \nabla \hat{\yb}^h)dx^h = \frac{1}{2}\int_\Omega Q_3 (x,\Gbb(x))dx
\end{equation}

\end{theorem}
\begin{prooof}
Let us first assume that $u,w,v_k, \etab$ are smooth. Then we
define for $(x_1,x_2^h,x_3^h) \in \bar{\Omega}^h$:
\begin{eqnarray} \nonumber
 \hat{\yb}^h (\Thetab^h(x_1,x_2^h,x_3^h)) &=& 
\Thetab^h (x_1,x_2^h,x_3^h)+\left(\begin{array}{c} h^2 u (x_1) \\ h v_2(x_1)
\\ h v_3(x_1)
\end{array} \right)\\ \nonumber  & &+h^2 \left(\begin{array}{c} -x_2 (v_2'p_2+v_3'p_3)(x_1)-x_3 (v_3'p_2-v_2'p_3)(x_1) \\ -x_2 (p_3 w)(x_1)-x_3 (p_2 w)(x_1) \\ x_2 (p_2w)(x_1)-x_3 (p_3w)(x_1) \end{array} \right)\\ & & + h^3\etab(x_1,\frac{x_2^h}{h}, \frac{x_3^h}{h}),
\end{eqnarray}
where $\etab: \Omega \to \ZR^3$ is going to be chosen later.
The convergence (a)-(d) and that $\etab^h \rightharpoonup
\etab$ weakly in $L^2(\Omega)$ and $\partial_k \etab^h
\rightharpoonup \partial_k \etab$ weakly in $L^2(\Omega)$ can
easily seen to be valid for this sequence. We also have
\begin{eqnarray} \nonumber
\nabla \hat{\yb}^h \nabla \Thetab^h &=& \nabla \Thetab^h + \left(\begin{array}{ccc} h^2 u' & -h(v_2'p_2+v_3'p_3) & -h(v_3'p_2-v_2'p_3)\\ hv_2' & -hp_3 w& -hp_2 w \\ hv_3' & hp_2 w &-hp_3 w\end{array} \right)\\ \nonumber  & &
+h^2 \left( \begin{array}{c} -x_2 (v_2'p_2+v_3'p_3)'+x_3(v_2'p_3-v_3'p_2)' \\ -x_2(p_3w)'-x_3 (p_2w)'\\
x_2(p_2w)'-x_3 (p_3w)'\end{array} \ \Bigg| \  \partial_2 \etab \ \Bigg| \ \partial_3 \etab \right)\\ \label{nablayt} & & +O(h^3).
\end{eqnarray}
From (\ref{nablayt}), by using (\ref{inverz}), we conclude
\begin{eqnarray} \nonumber
\nabla \hat{\yb}^h &=& \Ibb+h \left(\begin{array}{ccc} 0 & -v_2' & -v_3'\\ v_2' & 0& -w \\ v_3' &  w &0\end{array} \right)\\ & & \nonumber+h^2 \left(\begin{array}{ccc} u'+v_2'\theta_2'+v_3'\theta_3' & 0 & 0 \\ w \theta_3'  & v_2' \theta_2'& v_2'\theta_3' \\ -w \theta_2'& v_3'\theta_2' & v_3'\theta_3'   \end{array} \right)
\\ & & \nonumber+h^2 \left(\begin{array}{c}-x_2 (v_2''p_2+v_3''p_3)+x_3(v_2''p_3-v_3''p_2) \\ -x_2(p_3w')-x_3 (p_2w')\\
x_2(p_2w')-x_3 (p_3w')\end{array} \ \Bigg| \ \partial_2 \kab \ \Bigg| \  \partial_3 \kab \right)\\ & &  +O(h^3).
\end{eqnarray}
Using the identity $(\Ibb+\Mbb)^T (\Ibb+\Mbb)=\Ibb+2 \sym \Mbb
+\Mbb^T \Mbb$ we obtain
\begin{eqnarray} \nonumber
(\nabla \hat{\yb}^h)^T (\nabla \hat{\yb}^h) &=& \Ibb+2h^2 \sym \Jbb+2h^2 \sym \Kbb+h^2 \Abb^T \Abb+O(h^3),
\end{eqnarray}
where $\|O(h^3)\|_{L^\infty(\Omega;\ZR^{3 \times 3})} \leq
Ch^3$, for some $C>0$.

Taking the square root we obtain
\begin{equation} \label{strainzaul}
[(\nabla \hat{\yb}^h)^T (\nabla \hat{\yb}^h)]^{1/2}= \Ibb+h^2\widetilde{\Gbb}+O(h^3).
\end{equation}
We have $\det (\nabla \hat{\yb}^h)>0$ for sufficiently small
$h$. Hence by frame-indifference $W(x,(\nabla \hat{\yb}^h)\circ
\Thetab^h\circ P^h)= W(x, [\nabla \hat{\yb}^h)^T (\nabla
\hat{\yb}^h)]^{1/2}\circ \Thetab^h \circ P^h)$; thus by
(\ref{strainzaul}) and Taylor expansion we obtain:
$$ \frac{1}{h^4} W(x,(\nabla \hat{\yb}^h) \circ \Thetab^h \circ P^h ) \to\frac{1}{2} Q_3 (x,\widetilde{\Gbb}(x))\ \textrm{a.e}
$$
and by the property ii) of $W$ for $h$ small enough
$$ \frac{1}{h^4} W(x,(\nabla \hat{\yb}^h)\circ\Thetab^h \circ P^h) \leq \frac{1}{2}
C  (\| \Jbb \|^2+\| \Kbb \|^2+\| \Abb \|^4)+Ch. $$ The equality (\ref{upperjed}) follows by the
dominated convergence theorem. Namely, we have
\begin{eqnarray*}& & \frac{1}{h^6} \int_{\hat{\Omega}^h} W^h(x,
\nabla \hat{\yb}^h)dx=\frac{1}{h^4}\int_{\Omega} W(x,(\nabla
\hat{\yb}^h) \circ \Thetab^h \circ P^h)dx\\& & \hspace{25ex}
\to \frac{1}{2} \int_{\Omega} Q_3(x,\widetilde{\Gbb})
dx.\end{eqnarray*} In the general case, it is enough to
smoothly approximate $u,w$ in the strong topology of $W^{1,2}$,
$v_k$ in the strong topology of $W^{2,2}$, and $\etab,
\partial_k \etab$ in the strong topology of $L^2$ and to use
the continuity of the right hand side of (\ref{upperjed}) with
respect to these convergences.
\end{prooof}
\begin{remark} \label{zaidentifikaciju}
Notice that
\begin{eqnarray*}
\Kbb=\left( \Abb'\left(\begin{array}{c} 0 \\ x_2' \\ x_3'  \end{array}\right)\ | \ \partial_2 \kab \ | \ \partial_3  \kab \right)=\Lbb+\left( \Abb'\left(\begin{array}{c} 0 \\ x_2' \\ x_3'  \end{array}\right)\ | \ \partial_2 \betab \ | \ \partial_3  \betab \right).
\end{eqnarray*}
Here $\betab=\gamb \circ (\xb')^{-1}$ and
\begin{equation}
\Lbb=\left(\begin{array}{ccc} 0 & 0 & 0\\ 0& g_2^2 & g_2^3 \\ 0 & g_3^2 & g_3^3
\end{array} \right).
\end{equation}
From the fact that $\gamb \in \mathcal{C}$ we can conclude
$\betab \in \mathcal{B}$, where
\begin{eqnarray} \nonumber
\mathcal{B}&=&\{ \betab \in L^2(\Omega';\ZR^3): \int_\omega \betab=0, \ \partial_2 \betab, \partial_3 \betab \in L^2(\Omega'; \ZR^3), \\
& & \hspace{-4ex} \int_{\omega'(x_1)} (x_3' \betab_2(x_1,\cdot)-x_2' \betab_3(x_1,\cdot)dx'_2 dx'_3=0, \ \textrm{for a.e. } x_1 \in (0,L) \}.
\end{eqnarray}
\end{remark}
\subsection{Identification of the $\Gamma$-limit}
Let $Q:(0,L) \times \ZR \times \so(3) \to [0, +\infty)$ be
defined as
\begin{eqnarray} \nonumber
& &Q(x_1,t,\Fbb) = \\ & & \nonumber \min_{\alb \in W^{1,2}(\omega'(x_1);\ZR^3 )}\int_{\omega'(x_1)} Q_3\left( x, \left(\Fbb \left(\begin{array}{c}0 \\ x_2' \\ x_3' \end{array}\right)+t \eb_1 \Bigg| \partial_2 \alb \Bigg| \partial_3 \alb\right)\right)dx_2' dx_3',\\ & &  \label{minimum}
\end{eqnarray}
where $Q_3$ is the quadratic form defined in (\ref{defQ3}). For
$u, w \in W^{1,2}(0,L)$ and $v_2,v_3 \in W^{2,2}(0,L)$ we
introduce the functional
\begin{equation} \label{defi0}
I^0(u,v_2,v_3,w):=\frac{1}{2} \int_0^L Q(x_1, u'+v_2' \theta_2'+v_3'\theta_3'+\frac{1}{2} ((v_2')^2+(v_3')^2),\partial_1 \Abb)dx_1,
\end{equation}
where $\Abb \in W^{1,2}((0,L);\so(3))$ is defined by
(\ref{defAbb}). We shall state the result of
$\Gamma$-convergence of the functionals $\frac{1}{h^4} I^h$ to
$I^0$. Before stating the theorem we analyze some properties of
the limit density $Q$.
\begin{remark} \label{moraremark}
By using the remarks in the beginning of chapter 4 in
\cite{Mora1} the following facts can be concluded:
\begin{enumerate}[a)]
\item The functional $Q_3(x,\Gbb)$ is  coercive on
    symmetric matrices i.e. there exists a constant $C>0$,
    independent of $x$, such that $Q_3(x,\Gbb)\geq C\| \sym
    \Gbb\|^2$, for every $\Gbb$ (this is the direct
    consequence of the assumption iv) on $W$). The minimum
    in (\ref{minimum}) is attained. Since the functional
    $Q_3(x,\Gbb)$ depends only on the symmetric part of
    $\Gbb$, it is invariant under transformation $\alb
    \mapsto \alb+c_1+c_2(x')^\bot$ and hence the minimum
    can be computed on the subspace
    \begin{eqnarray*} V_{x_1}&:=&\Bigg\{ \alb \in W^{1,2}
    (\omega'(x_1), \ZR^3): \ \int_{\omega'(x_1)} \alb=0,\\ & & \hspace{8ex}
    \int_{\omega'(x_1)} (x_3' \alb_2- x_2' \alb_3)dx_2'
    dx_3'=0 \Bigg\}.
    \end{eqnarray*}
Strict convexity of $Q_3(x, \cdot)$ on symmetric matrices
ensures that the minimizer is unique in $V$.
\item Fix $x_1 \in (0,L)$, $t \in \ZR$ and $\Fbb \in
    \so(3)$. Let $\alb^{min} \in V$ be the unique minimizer
    of the problem (\ref{minimum}). We set
    $$g(x_2',x_3')= \Fbb \left(\begin{array}{c}0 \\ x_2' \\ x_3' \end{array}\right)+t \eb_1, \quad b_{ij}^{hk}=\frac{\partial^2 W}{\partial \Fbb_{ih} \partial \Fbb_{jk}} (x, \Ibb), $$
and we call $B^{hk}$ the matrix in $\ZR^{3 \times 3}$ whose
elements are given by $(B^{hk})_{ij}=b^{hk}_{ij}$. Then
$\alb^{min}$ satisfies the following Euler-Lagrange
equation:
\begin{equation}
\int_{\omega'(x_1)} \sum_{h,k=2,3} (B^{hk} \partial_k \alb^{min}, \partial_h \varphi)dx_2'dx_3'=-\int_{\omega'(x_1)} \sum_{h=2,3} (B^{h1}g, \partial_h \varphi) dx_2'dx_3',
\end{equation}
for every $\varphi \in W^{1,2} (\omega'(x_1);\ZR^{3 \times
3})$.  From this equation it is clear that $\alb^{min}$
depends linearly on $(t,\Fbb)$. Moreover $Q$ is uniformly
positive definite, i.e.
\begin{equation}
Q(x_1,t, \Fbb) \geq C (t^2+\| \Fbb \|^2), \quad \forall t \in \ZR, \forall \Fbb \in \so(3),
\end{equation}
and the constant $C$ does not depend on $x_1$.
\item By mimicking the proof of Remark 4.3 in \cite{Mora1}
    it can be seen that there exists a constant $C'$
    (independent of $x_1$, $t$ and $\Fbb$) such that
   \begin{equation}
\| \partial_2 \alb^{min} \|_{L^2(\omega'(x_1);\ZR^{3 \times 3})}+\| \partial_3 \alb^{min} \|_{L^2(\omega'(x_1);\ZR^{3 \times 3})} \leq C'\|g\|^2_{L^2(\omega'(x_1); \ZR^{3 \times 3})},
    \end{equation}
for a.e. $x_1 \in (0,L)$. To adapt the proof we only need
to have that the constant in the Korn's inequality
\begin{equation}
\int_{\omega'(x_1)} \sum_{j,k=2,3} | \partial_k \alb^{min}_j |^2dx_2' dx_3' \leq C_1 \int_{\omega'(x_1)} \sum_{j,k=2,3} |e_{jk} (\alb^{min}) |^2 dx_2' dx_3'
\end{equation}
can be chosen independently of $x_1$. This is proved in
Lemma \ref{Kooornova}.
\item When $Q_3$ does not depend on $x_2$, $x_3$  we can
    find a more explicit representation for $Q$. More
    precisely $Q$ can be decomposed into the sum of two
    quadratic forms
    $$ Q(x_1,t,\Fbb)=Q_1(x_1,t)+Q_2(x_1,\Fbb),$$
where
\begin{eqnarray}
Q_1(x_1,t)&:=& \min_{\ab,\bb \in \ZR^3} Q_3(x_1,(t \eb_1 | \ab | \bb)),\\
Q_2(x_1,0,\Fbb) &:=& Q(x_1,0,\Fbb).
\end{eqnarray}
The relations (\ref{nulaaa}) are only needed for this. If
we assume the isotropic and homogenous case i.e.
$$ Q_3(\Fbb)=2 \mu \Big| \frac{\Fbb+\Fbb^T}{2} \Big|^2+\lambda
(\tra \Fbb)^2, $$ then after some calculation (see Remark 3.5 in \cite{Mora0}) it can be shown that
\begin{eqnarray*}
Q_1 (t)&=& \frac{\mu (3 \lambda+2 \mu)}{\lambda+\mu}t^2 \\
Q_2(x_1,\Fbb) &=& \frac{\mu (3 \lambda+2 \mu)}{\lambda+\mu} \Big(\Fbb_{12} \int_{\omega'(x_1)} (x_2')^2 dx_2'dx_3'\\ & & +2\Fbb_{12} \Fbb_{13} \int_{\omega'(x_1)} x_2'x_3' dx_2'dx_3'+\Fbb_{13} \int_{\omega'(x_1)} (x_3')^2 dx_2'dx_3'   \Big)\\ & & +\mu \tau  \Fbb_{23},
\end{eqnarray*}
where the constant $\tau$ is so-called torsional rigidity,
defined as
$$ \tau (\omega'(x_1))= \tau(\omega) =\int_{\omega}
(x_2^2+x_3^2-x_2 \partial_3 \varphi+x_3 \partial_2 \varphi)
dx_2 dx_3,$$
and $\varphi$ is the torsion function i.e. the
solution of the Neumann problem $$\left\{ \begin{array}{ll}
\nabla \varphi =0 &\  \textrm{in } \omega \\ \partial_{\nu}
\varphi=-(x_3,-x_2) \cdot \nu & \ \textrm{on } \partial
\omega
\end{array} \right. $$
\end{enumerate}
\end{remark}
The following theorem can be proved in the same way as Theorem
4.5 in \cite{Mora} (we need Theorem \ref{prepisano}, Lemma
\ref{identifikacija}, Theorem \ref{upperbound}, Remark
\ref{zaidentifikaciju} and Remark \ref{moraremark}).
\begin{theorem} \label{najglavnijiii}
As $h \to 0$, the functionals $\frac{1}{h^4} I^h$ are
$\Gamma$-convergent to the functional $I^0$ given in
(\ref{defi0}), in the following sense:
\begin{enumerate}[i)]
\item (compactness and liminf inequality) if \  $\limsup_{h
    \to 0} h^{-4} I^h< +\infty$ then there exists constants
    $\oRbb^h \in \SO(3)$, $c^h \in \ZR^3$ such that (up to
    subsequences) $\oRbb^h \to \oRbb$ and the functions
    defined by
\begin{eqnarray*}
& & \widetilde{\yb}^h:= (\oRbb^h)^T\yb^h-c^h, \quad u^h = \frac{1}{A}\int_\omega \frac{\widetilde{\yb}_1^h\circ \Thetab^h\circ P^h-x_1}{h^2}dx_2 dx_3 \\
& & v^h_k = \frac{1}{A}\int_\omega \frac{\widetilde{\yb}_k^h\circ \Thetab^h\circ P^h-h\theta_k}{h} dx_2 dx_3 \\
& & w^h = \frac{1}{A\mu(\omega)}\int_\omega \frac{x_2'(\widetilde{\yb}_3\circ \Thetab^h\circ P^h)-x_3'(\widetilde{\yb}_2\circ \Thetab^h\circ P^h)}{h^2} dx_2 dx_3
\end{eqnarray*}
satisfy
\begin{enumerate}
\item  $(\nabla \widetilde{\yb}^h) \circ \Thetab^h
    \circ P^h \to \Ibb$ in $L^2(\Omega)$.
\item there exist $u,w \in W^{1,2}(0,L)$ such that $u^h
    \rightharpoonup u$ and $w^h \rightharpoonup w$
    weakly in $W^{1,2}(0,L)$.
\item there exists $v_k \in W^{2,2}(0,L)$ such that
    $v_k^h \to v_k$ strongly in $W^{1,2}(0,L)$ for
    $k=2,3$.
\end{enumerate}
Moreover we have
\begin{equation}
\liminf_{h \to 0} \frac{1}{h^4} I^h(\yb^h)\geq I^0(u,v_2,v_3,w).
\end{equation}
\item (limsup inequality) for every $v,w \in W^{1,2}(0,L)$,
    $v_2,v_3 \in W^{2,2}(0,L)$ there exists $(\hat{\yb}^h)$
    such that (a)-(c) hold (with $\widetilde{\yb}^h$
    replaced by $\hat{\yb}^h$) and
    \begin{equation}
\lim_{h \to 0} \frac{1}{h^4} I^h(\hat{\yb}^h) = I^0 (u,v_2,v_3,w)
    \end{equation}
\end{enumerate}
\end{theorem}
\begin{remark} \label{napcijela1}
Let $f_2, f_3 \in L^2(0,L)$. We introduce the functional
\begin{equation}
J^0=I^0(u,v_2,v_3,w)-\int_0^L \sum_{k=2,3} f_k v_k,
\end{equation}
for every $u \in W^{1,2}(0,L)$, $v_2, v_3 \in W^{2,2}(0,L)$,
and $w \in W^{1,2}(0,L)$. The functional $J^0$ can be obtained
as $\Gamma$-limit of the energies $\frac{1}{h^4} I^h$ by adding
a term describing transversal body forces of order $h^3$ (see
\cite{Muller3}, see also \cite{Velcic}). For longitudinal body
forces see \cite{Muller5}. The problem for longitudinal body
forces arises because the longitudinal forces should be of
order $h^2$, the same order as for the model in \cite{Mora0}.
One needs to impose certain stability condition to see which
model describes the behavior of the body for the longitudinal
forces of order $h^2$.
\end{remark}
\begin{remark}
The term $u'+v_2'\theta_2'+v_3' \theta_3'+\frac{1}{2}
((v_2')^2+(v_3')^2)$ in the strain measures the extension of
the central line (which is of the second order). Namely, if we
approximate the deformation of the weakly curved rod by:
\begin{eqnarray}
\varphi_1 (x_1,x_2,x_3)&=& x_1+h^2u+h^2 x_2'(v_2'+\theta_2')+h^2x_3'(v_3'+\theta_3') \\
\varphi_k (x_1,x_2,x_3)&=& h\theta_k+hx_k'+hv_k+h^2 (x_k')^\bot w, \ \textrm{for } k=2,3,
\end{eqnarray}
we see, that it is valid
\begin{eqnarray*}
\| \partial_1 \varphib (x_1,0,0)\|^2-\|\partial_1 \Thetab^h(x_1,0,0)\|^2&=&h^2\Big(2u'+2v_2'\theta_2'+2v_3'\theta_3'\\ & &+(v_2')^2+(v_3')^2\Big).
\end{eqnarray*}
\end{remark}
\begin{remark} \label{napcijela2}
The existence of the solution for the functional $J^0$ under
the Dirichlet boundary condition for $v_k$ at both ends of the
rod can be proved directly. It is also enough that we impose
$v_2,v_2',v_3,v_3'$ at the one end. The existence can also be
proved for the free boundary condition under the hypothesis
that $\int_0^L \fb_k dx_1=0$, $\int_0^L x_1 \fb_k dx_1=0$ for
$k=2,3$. It can be done in the same way as the proof of Lemma 5
in \cite{Velcic}.
\end{remark}

\end{document}